\documentclass [12pt,a4paper,reqno]{amsart}
 \input amssymb.sty
\newcommand{\ef}{\end{equation}}
\chardef\bslash=`\\ % p. 424, TeXbook
\newtheorem{thm}{Theorem} [section]
\newtheorem*{thm*}{Theorem}
\newtheorem{cor}[thm]{Corollary}
\newtheorem{lem}[thm]{Lemma}

 \newtheorem{defn}[thm]{Definition}

\newtheorem{remark}[thm]{Remark}
\newtheorem{acknowledgment*}[thm] {Acknowledgment}

\newtheorem*{acknowledgment} {Acknowledgment}
\newcommand{\thmref}[1]{Theorem~\ref{#1}}

\newcommand{\lemref}[1]{Lemma~\ref{#1}}
\newcommand{\corref}[1]{Corollary~\ref{#1}}

 \renewcommand{\sectionmark}[1]{}

\newcommand{\iy}{\infty}

\newcommand{\doe}{\overset{\text{def}}{=}}
 
\newcommand{\loc} {\operatorname{loc}}

 \date{}
\begin{document}

\title[Davies-Harrell Representations, Otelbaev's Inequalities...]{Davies-Harrell Representations, Otelbaev's Inequalities\\
and properties of solutions of Riccati equations}
\author[N.A. Chernyavskaya]{N.A. Chernyavskaya}\address{
Department of Mathematics and Computer Science, Ben-Gurion University of the
Negev, P.O.B. 653, Beer-Sheva, 84105, Israel}
\author[L.A. Shuster]{L.A. Shuster}
 \address{Department of Mathematics,
 Bar-Ilan University, 52900 Ramat Gan, Israel}

\begin{abstract}
We consider  an equation
\begin{equation}\label{1}
y''(x)=q(x)y(x),\quad x\in R
\end{equation}
under the following assumptions on $q(x):$
\begin{equation}\label{2}
0\le q(x)\in L_1^{\loc}(R),\quad \int_{-\infty}^xq(t)dt>0,\quad \int_x^\infty
q(t)dt>0\quad\text{for all}\quad x\in R.
\end{equation}
Let $v(x)$ (resp. $u(x)$) be a positive non-decreasing (resp. non-increasing)
solution of \eqref{1} such that
$$v'(x)u(x)-u'(x)v(x)=1,\quad x\in R.$$
These properties determine $u(x)$ and $v(x)$ up to mutually inverse positive
constant factors, and the function $\rho(x)=u(x)v(x),$\ $x\in R$ is uniquely
determined by $q(x)$,\ $x\in R.$ In the present paper, we obtain an asymptotic
formula for computing $\rho(x)$ as $|x|\to\infty.$ As an application, under
conditions \eqref{2}, we study the behavior at infinity of solution of the
Riccati equation
$$ z'(x)+z(x)^2=q(x),\quad x\in R.$$
\end{abstract}
\maketitle

\baselineskip 20pt

\section{Introduction}\label{introduction}
\setcounter{equation}{0} \numberwithin{equation}{section}

In the present paper, we consider an equation
\begin{equation}\label{1.1}
y''(x)=q(x)y(x),\quad x\in R
\end{equation}
under assumptions
\begin{equation}\label{1.2}
0\le q(x)\in L_1^{\loc}(R),\quad \int_{-\infty}^xq(t)dt>0,\quad \int_x^\infty
q(t)dt>0\quad\text{for all}\quad x\in R.
\end{equation}
Further, we assume conditions \eqref{1.2} are satisfied, without special
mention.  Our general goal is to study some asymptotic properties (as
$|x|\to\infty$) of solution of equations \eqref{1.1}.

In order to give a concrete statement of the problem, we need the following
known facts (see \cite[Ch.\,XI, \S6]{1},\cite{2}). First we note that equation
\eqref{1.1} has a fundamental system of solutions (FSS) $\{u(x),v(x)\}$ which
is defined, up to mutually inverse positive constant factors, by the properties
\begin{equation}\label{1.3}
v(x)>0,\quad u(x)>0,\quad v'(x)\ge 0,\quad u'(x)\le 0\quad\text{for}\quad x\in
R,
\end{equation}
\begin{equation}\label{1.4}
v'(x)u(x)-u'(x)v(x)=1\quad\text{for}\quad x\in R,
\end{equation}
\begin{equation}\label{1.5}
\lim_{x\to-\infty}\frac{v(x)}{u(x)}=\lim_{x\to\infty}\frac{u(x)}{v(x)}=0,
\end{equation}
\begin{equation}\label{1.6}
\int_{-\infty}^0\frac{d\xi}{v(\xi)^2}=\int_0^\infty\frac{d\xi}{u(\xi)^2}=\infty,\quad
\int_{0}^\infty\frac{d\xi}{v(\xi)^2}<\infty,\quad
\int_{-\infty}^0\frac{d\xi}{u^2(\xi)}<\infty.
\end{equation}
Relations \eqref{1.3}--\eqref{1.6} mean that $u(x)$ and $v(x)$ are principal
solutions of \eqref{1.1} on $(0,\infty)$ and $(-\infty,0),$ respectively (see
\cite[Ch.\,XI, \S6]{1}. Therefore, we call an FSS $\{(u(x),v(x)\}$ with
properties \eqref{1.3}--\eqref{1.6} a principal FSS (PFSS) of equation
\eqref{1.1}.

The solutions $u(x),v(x)$ from a PFSS of \eqref{1.1} are related as follows:
\begin{equation}\label{1.7}
u(x)=v(x)\int_{x}^\infty\frac{dt}{v(t)^2},\quad
v(x)=u(x)\int_{-\infty}^x\frac{dt}{u(t)^2},\quad x\in R.
\end{equation}
{}From \eqref{1.7}, it follows that the function $\rho(x)$
\begin{equation}\label{1.8}
\rho(x)\doe
u(x)v(x)=v(x)^2\int_x^\infty\frac{dt}{v(t)^2}=u(x)^2\int_{-\infty}^x\frac{dt}{u(t)^2},\quad
x\in R
\end{equation}
does not depend on the choice of a PFSS and is determined uniquely by equation
\eqref{1.1}, i.e., by the function $q(x).$ Therefore, the Davies-Harrell
representation \eqref{1.9} for a PFSS of equation \eqref{1.1} (see \cite{3}) is
very important for the theory of equation \eqref{1.1}
\begin{equation}\label{1.9}
u(x)=\sqrt{\rho(x)}\exp\left(-\frac{1}{2}\int_{x_0}^x\frac{d\xi}{\rho(\xi}\right),\quad
v(x)=\sqrt{\rho(x)}\exp\left(\frac{1}{2}\int_{x_0}^x\frac{d\xi}{\rho(x)}\right),\quad
x\in R.
\end{equation}
Here $x_0$ is the unique root of the equation $u(x)=v(x)$ (such an
interpretation of Davies-Harrell's formulas was proposed in \cite{2}). Thus for
all $x\in R$, any  PFSS of \eqref{1.1} can be expressed via $\rho(x)$, and the
choice of a particular PFSS of \eqref{1.1} is determined by the choice of $x_0$
in \eqref{1.9}.

Representation \eqref{1.9} becomes even more important if one takes into
account Otelbaev's a~priori inequalities
\begin{equation}\label{1.10}
\frac{d(x)}{4}\le \rho(x)\le\frac{3}{2}d(x),\quad x\in R.
\end{equation}
Here $d(x)$ is the unique solution in $d\ge0$ of
\begin{equation}\label{1.11}
d\int_{x-d}^{x+d}q(t)dt=2,\quad x\in R.
\end{equation}

\begin{remark}\label{rem1.1} {\rm The function $d(x)$ was introduced by M.
Otelbaev (see, for example, \cite{7}). It is well-defined (see \S2,
\lemref{lem2.1}). Inequalities of type \eqref{1.10} were first obtained in
\cite{4} (under requirements of $q(x)$ stronger than \eqref{1.2}), and
therefore we relate them and the function $d(x)$ to M. Otelbaev. Note that in
\cite{4} another auxiliary function, more complicated than $d(x)$, was used.
See \cite{2} for the proof of estimates \eqref{1.10} under conditions
\eqref{1.2}.  }
\end{remark}

The study of $\rho(x)$ started in \cite{2} was continued in \cite{5,6}. In
\cite{5} more precise inequalities of type \eqref{1.10} were obtained, and in
\cite{6}, under some additional requirements of \eqref{1.2} to $q(x),$ an
asymptotic formula for computation of $\rho(x)$ as $|x|\to\infty,$ was
obtained. We shall need this formula later. To state it, let us introduce the
following

\begin{defn}\label{defn1.2}\cite{6}
{\rm Suppose that condition \eqref{1.2} holds. We say that $q(x)$ belongs to
the class $H$ (and write $q(x)\in H$) if there exists a continuous function $k
(x)$ in $x\in R$ with properties:}
 \begin{enumerate} \item[1)]
\begin{equation}\label{1.12}  k (x)\ge 2,\quad x\in R;\ k (x)\to\infty\quad\text{as}\quad
|x|\to\infty; \end{equation} \item[2)] {\rm there is an absolute positive
constant $c_1$ such that for all $x\in R$ the following inequalities hold:}
\begin{equation}\label{1.13}
c_1^{-1}k (x)\le k (t)\le c_1k (x)\quad\text{for}\quad t\in [x-k (x)d(x),x+k
(x)d(x)]
\end{equation}
\item[3)] {\rm there is an absolute positive constant $c_2$ such that for all
$x\in R$ the following estimate holds:}
\begin{equation}\label{1.14}
\Phi(x)\doe k (x)d(x)\sup\limits_{z\in[0,k
(x)d(x)]}\left|\int_0^z(q(x+t)-q(x-t))dt \right| \le c_2.
\end{equation}
\end{enumerate}
\end{defn}

In the sequel we assume that if $q(x)\in H,$ we denote by $k (x)$ the function
from Definition \ref{defn1.2}. For example, if $q(x)\in H,$ then  below by
$F(x)$ we denote the function
\begin{equation}\label{1.15}
F(x)=\doe \sqrt{k (x)}d(x)\sup_{z\in[0,\sqrt{k
(x)}d(x)]}\left|\int_0^z(q(x+t)-q(x-t))dt\right|,\quad x\in R.
\end{equation}
Later we omit the reference to the conditions \eqref{1.12}, \eqref{1.13} and
\eqref{1.14},   which the function $k (x)$ in \eqref{1.15} (and  in any similar
situation) satisfies. By $c$ we denote any absolute positive constants which
are not essential for exposition and which may differ even within a single
chain of computations. Constants essential for exposition are supplied with
indices, as, for example, in Definition \ref{defn1.2}.

\begin{thm}\label{thm1.3} \cite{6}
Suppose that $q(x)\in H$ and, in addition, $q(x)\ge 1$ for $x\in R.$ Then for
all $|x|\gg 1,$ we have
\begin{equation}\label{1.16}
|\rho'(x)|\le c[F(x)+\exp(c^{-1}\sqrt{k (x)})]\le \frac{c}{\sqrt{k (x)}},
\end{equation}
\begin{equation}\label{1.17}
\rho(x)=\frac{d(x)}{2}(1+\varepsilon(x)),\quad |\varepsilon(x)|\le
c\alpha(x)\le\frac{c}{\sqrt{k (x)}}.
\end{equation}
Here (see \eqref{1.15}):
\begin{equation}\label{1.18}
\alpha(x)=\begin{cases} \exp(-c^{-1}\sqrt{k (x)})+\sup\limits_{t\ge
x-d(x)} F(t),\quad& \text{if}\quad x\ge 0\\
\exp(-c^{-1}\sqrt{k (x)})+\sup\limits_{t\le x+d(x)}F(t),\quad&\text{if}\quad
x\le 0.
\end{cases}
\end{equation}
\end{thm}

Here are some comments on \thmref{thm1.3}. The main goal of this statement is
to make inequalities \eqref{1.10} more precise for $|x|\gg 1.$ A solution is
suggested in \eqref{1.17}--\eqref{1.18}. Clearly, in view of representations
\eqref{1.19} of PFSS, formulas of type \eqref{1.17}--\eqref{1.18} are important
for the theory of equation \eqref{1.1}. In addition, they are applied, for
example, in the spectral theory of the Sturm-Liouville operator and in the
theory of the Riccati equation (see \cite{5,6}). Therefore, their further
development may be useful for equation \eqref{1.1} as well as for its
applications. Note that relations \eqref{1.17}--\eqref{1.18} and
\eqref{1.10}--\eqref{1.11} do not completely agree with one another. In
particular, Otelbaev's inequalities are local because to estimate the function
$\rho(x)$ in a point $x\in R$, one only uses the values $q(t)$ for all $t$ from
the finite segment $[x-d(x),x+d(x)].$ In contrast, asymptotic estimates
\eqref{1.17} are not local because to estimate $\rho(\cdot)$ in a point $x$\
($|x|\gg1)$ one uses the values $q(t)$ for all $t$ belonging to one of the
infinite intervals $(-\infty,x+d(x)]$ or $[x-d(x),\infty)$ (see \eqref{1.11}
and \eqref{1.18}). Analysis of the examples to \thmref{thm1.3} from \cite{5,6}
shows that in formula \eqref{1.17}, the estimates of the remainder term
$\varepsilon(\cdot)$ in a point $x$\ $(|x|\gg1)$ are always formed from the
values $q(t)$ related to some local neighborhood of $x.$ This means that in
\eqref{1.17}--\eqref{1.18}, when estimating $\varepsilon(x),$ perhaps we impose
redundant conditions on the function $q(\cdot).$

Now, that we have clarified some disadvantages of the relations
\eqref{1.17}--\eqref{1.18}, we are able finally to formulate the main goal of
this paper:  to obtain an asymptotic formula with a local estimate of the
remainder term for computing $\rho(x)$ as $|x|\to\infty.$ This problem is
solved in \thmref{thm1.4} which is the main result of the present paper, as
follows:

\begin{thm}\label{thm1.4}
Suppose $q(x)\in H.$ Then for all $|x|\gg1$ estimates \eqref{1.16} hold, and we
have the following relations:
\begin{equation}\label{1.19}
\rho(x)=\frac{d(x)}{2}(1+\varepsilon(x)),\quad |\varepsilon(x)|\le c\beta(x).
\end{equation}
 Here
 \begin{equation}\label{1.20}
\beta(x)=\exp(-c^{-1}\sqrt{k (x)})+ \sup_{t\in\Delta(x)}F(x)\le\frac{c}{\sqrt{k
(x)}},\quad
 \Delta(x)=[x-d(x),x+d(x)].
 \end{equation}
\end{thm}

We give here \thmref{thm1.5} containing a more detailed variant of formula
\eqref{1.19} which is intended for the following particular application. We
plan  to apply \eqref{1.19} for constructing approximations to the solutions of
the equation
\begin{equation}\label{1.21}
-y''(x)+q(x)y(x)=f(x),\quad x\in R
\end{equation}
with $f(x)\in L_p(R),$\ $p\in[1,\infty]$\ $(L_\infty(R):=c(R))$ and $q(x)\in
H.$ To solve this problem, we need a more detailed version of \thmref{thm1.4}.
It is convenient to state it here as a separate assertion. First note that the
functions $q(x)\in H$ have the following property.  For every $q(x)\in H$ there
is an absolute positive constant $c_3$ such that (see \S2, \lemref{lem2.7})
\begin{equation}\label{1.22}
c_3^{-1}d(x)\le d(t)\le c_3d(x)\quad\text{for}\quad |t-x|\le\sqrt{k
(x)}d(x),\quad x\in R.
\end{equation}
We introduce some notation:
\begin{equation}\label{1.23}
d_0=\sup_{x\in R}d(x),
\end{equation}
\begin{equation}\label{1.24}
\eta_1(x)=4[F(x)+\sqrt{c_3}\exp(-(3c_3)^{-1}\sqrt{k (x)})],\quad x\in R,
\end{equation}
\begin{equation}\label{1.25}
\eta_2(x)=65\left[\sup_{t\in\Delta(x)}F(t)+\sqrt{c_3}\exp\left(-\left(3c_3\sqrt{c_3}\right)^{-1}\sqrt{k
(x)}\right)\right],\quad x\in R.
\end{equation}

\begin{thm}\label{thm1.5}
Suppose $q(x)\in H$ and $d_0<\infty.$ Denote by $S_0$ a point on the number
axis such that for all $|x|\ge s_0$ the following inequalities hold (see
\eqref{1.12} and \eqref{1.16}):
\begin{equation}\label{1.26}
k (x)\ge 64c_2^2,\quad \eta_1(x)\le 10^{-3},
\end{equation}
and set $s_1=s_0+d_0+1.$ Then the following relations hold:
\begin{equation}\label{1.27}
|\rho'(x)|\le\eta_1(x)\quad\text{for}\quad |x|\ge s_1,
\end{equation}
\begin{equation}\label{1.28}
\rho(x)=\frac{d(x)}{2}(1+\varepsilon(x)),\quad |\varepsilon(x)|\le
\eta_2(x)\quad \text{for}\quad |x|\ge s_1.
\end{equation}
\end{thm}

\begin{remark}\label{rem1.6} {\rm
Theorems \ref{thm1.4} and \ref{thm1.5} will be proved together in \S\S2--4.
Each section contains a separate part of the proof accompanied by necessary
comments. The proof of formula \eqref{1.19} for a concrete equation, along with
technical details, is contained in \S8.}
\end{remark}

Let us now compare Theorems \ref{thm1.3} and \ref{thm1.4}. \thmref{thm1.3}
contains the requirement $q(x)\ge1$ for $x\in R$ which is not contained in
\thmref{thm1.4}. This restriction is not essential since relations
\eqref{1.17}--\eqref{1.18} can also be obtained by the method of \cite{6} using
only condition \eqref{1.2} and the condition $q(x)\in H.$ Nevertheless, in
order to reveal the principal difference between Theorems \ref{thm1.3} and
\ref{thm1.4}, in our comments below we assume that the condition $q(x)\ge1,$\
$x\in R$ holds true. This convention immediately implies that a new ``local"
version of \eqref{1.19}--\eqref{1.20} of asymptotic estimates at infinity for
the function $\rho(x)$ is obtained under the same assumptions on which, in
\thmref{thm1.3} only guaranteed a ``non-local" form of estimates
\eqref{1.17}--\eqref{1.18}. This ``qualitative" advantage of \thmref{thm1.4}
with respect to \thmref{thm1.3} is evidently important for solving theoretical
problems related to the properties of the function $\rho(x)$ at infinity.
However, when applied to concrete equations, a ``general quantitative
advantage" of relations \eqref{1.19}--\eqref{1.20} with respect to
\eqref{1.17}--\eqref{1.18} turns out to be more important.

This advantage can be expressed as follows: the asymptotic formula
\eqref{1.19}--\eqref{1.20} can be viewed as a refinement of the asymptotic
formula \eqref{1.17}--\eqref{1.18} in the class $H.$ To justify that, note that
Theorems \ref{thm1.3} and \ref{thm1.4} differ only in the functions $\alpha(x)$
and $\beta(x)$ which give an estimate of the same remainder terms
$\varepsilon(x)$ in the asymptotic formula
\begin{equation}\label{1.29}
\rho(x)=\frac{d(x)}{2}(1+\varepsilon(x)),\quad
\lim_{|x|\to\infty}\varepsilon(x)=0.
\end{equation}
Here both functions are constructed by the function $q(x),$\ $x\in R,$ are
continuous for $x\in R$ and satisfy the relations
\begin{equation}\label{1.30}
0<\beta(x)\le\alpha(x),\quad x\in R,\qquad
\lim_{|x|\to\infty}\alpha(x)=\lim_{|x|\to\infty}\beta(x)=0.
\end{equation}

Denote
\begin{equation}\label{1.31}
L=\sup_{q(x)\in H}\varlimsup\limits_{|x|\to\infty}\frac{\alpha(x)}{\beta(x)}.
\end{equation}
We say that the asymptotic formulas \eqref{1.17}--\eqref{1.18} and
\eqref{1.19}--\eqref{1.20} are equivalent in the class $H$ if $L<\infty.$ If
$L=\infty$, we say that the asymptotic formula \eqref{1.19}--\eqref{1.20} is a
refinement of the asymptotic formula \eqref{1.17}--\eqref{1.18} in the class
$H.$ With this terminology, the following assertion give the main relationship
between Theorems \ref{thm1.3} and \ref{thm1.4}.

\begin{thm}\label{thm1.7}
The asymptotic formula \eqref{1.19}--\eqref{1.20} is a refinement of the
asymptotic formula \eqref{1.17}--\eqref{1.18} in the class $H.$
\end{thm}

We give here  an example of an application of \thmref{thm1.4}. Consider a
Riccati equation
\begin{equation}\label{1.32}
y'(x)+y(x)^2=q(x),\quad x\in R.
\end{equation}

In \S6, we prove the following theorem which complements one of the results of
\cite{5}.
\begin{thm}\label{thm1.8}
Suppose $q(x)\in H.$ Then the following assertions hold:
\begin{enumerate}
\item[A)] There exists a unique solution $y_1(x)\ (y_2(x))$ of equation
\eqref{1.32} defined for all $x\in R$ and satisfying the equalities
 \begin{equation}
\begin{aligned}\label{1.33}
& \lim_{x\to-\infty}y_1(x)d(x)
 =\lim_{x\to\infty}y_1(x)d(x)=-1  \\
 &\left(\lim_{x\to-\infty}y_2(x)d(x)=\lim_{x\to\infty}y_2(x)d(x)=1\right).
\end{aligned}
 \end{equation}
 \item[B)] Let $y_+(x)$ be a solution of \eqref{1.32} defined on $[c,\infty)$
 for some $c.$ Then $y_+(x)\ne y_1(x)$ if and only if
 \begin{equation}\label{1.34}
 \lim_{x\to\infty}y_+(x)d(x)=1.
 \end{equation}
 \item[C)] Let $y_-(x)$ be a solution of \eqref{1.32} defined on $(-\infty,c]$
 for some $c.$ Then $y_-(x)\ne y_2(x)$ if and only if
 \begin{equation}\label{1.35}
 \lim_{x\to-\iy}y_-(x)d(x)=-1.
 \end{equation}
\end{enumerate}
\end{thm}
Note that an example of \thmref{thm1.8} is contained in \S8.

\begin{acknowledgment}{\rm
The authors thank Prof. Ya.~M. Goltser and Prof. Z.~S. Grinshpun for useful
discussions. }
\end{acknowledgment}

\section{Technical assertions}

In this section, we present some auxiliary assertions on the properties of the
function $d(x)$ (see \eqref{1.11}). Most of these lemmas were obtained in
\cite{6} under the assumption
\begin{equation}\label{2.1}
1\le q(x)\in L_1^{\loc}(R),\quad x\in R.
\end{equation}
To pass from condition \eqref{2.1} to condition \eqref{1.2}, we have to prove
that the ``old" assertions remain true under the ``new" assumptions. Our new
proofs are simpler and shorter than the previous ones and significantly differ
from those presented in \cite{6}.

\begin{lem}\label{lem2.1} \cite[Ch.I, \S5]{7} \ For every given $x\in R,$
equation \eqref{1.1} has a unique solution in $d\ge0.$
\end{lem}

\begin{proof} The functions
$$\varphi_1(d)=\frac{2}{d},\quad \varphi_2(d)=\int_{x-d}^{x+d}q(t)dt,\quad
d\in(0,\infty)$$ have the following properties:
\begin{enumerate}\item[1)] the function $\varphi_1(d)$ is monotone decreasing
from infinity to zero on $(0,\iy);$ \item[2)] the function $\varphi_2(d)$ is
non-decreasing and non-negative on $(0,\iy)$ and, in addition,
$\lim\limits_{d\to\iy}\varphi_2(d)=\iy$ (see \eqref{1.2}).
\end{enumerate}
{}From 1)--2) and the continuity of the two functions, it follows that their
graphs intersect at one point.
\end{proof}

\begin{lem}\label{lem2.2} \  For every  $x\in R,$
the inequality $\eta\ge d(x)$\ $(0\le\eta\le d(x))$ holds if and only if
\begin{equation}\label{2.2}
S(\eta)\ge 1\ (S(\eta)\le 2),\quad S(\eta)\doe\eta\int_{x-\eta}^{x+\eta}q(t)dt.
\end{equation}
\end{lem}

\renewcommand{\qedsymbol}{}
\begin{proof}[Proof of \lemref{lem2.2}] Necessity.

If $\eta\ge d(x),$ then $S(\eta)\ge S(d(x))=2.$
\end{proof}

\renewcommand{\qedsymbol}{\openbox}
\begin{proof}[Proof of \lemref{lem2.2}] Sufficiency.

Assume the contrary: $S(\eta)\ge 2$, but $\eta<d(x).$ Then $2\le S(\eta)\le
S(d(x))=2\ \Rightarrow\ S(\eta)=2.$ Hence $\eta=d(x)$ by \lemref{lem2.1}.
Contradiction.
\end{proof}

For a given $x\in R,$ consider an equation in $d\ge0:$
\begin{equation}\label{2.3}
G(d)=1,\qquad G(d)\doe\int_0^d\int_{x-t}^{x+t}q(\xi)d\xi dt,\qquad d\ge0.
\end{equation}

\begin{lem}\label{lem2.3}
For every $x\in R,$ equation \eqref{2.3} has a unique positive solution. Denote
it by $\hat d(x).$ The function $\hat d(x)$ satisfies the inequalities
\begin{equation}\label{2.4}
d(x)\le 2\hat d(x)\le 3d(x),\quad x\in R,
\end{equation}
\begin{equation}\label{2.5}
1\le\hat d(x)\int_{x-\hat d(x)}^{x+\hat d(x)}q(t)dt,\quad x\in R.
\end{equation}
In addition, $\hat d(x)$ has a continuous derivative for $x\in R,$ and
\begin{equation}\label{2.6}
|\hat d'(x)|\le \hat d(x)\left|\int_0^{\hat d(x)}(q(x+t)-q(x-t))dt\right|,\quad
x\in R.
\end{equation}
\end{lem}

 \begin{remark}\label{rem2.4} {\rm The function $\hat d(x)$ was introduced in
 \cite{8} under condition \eqref{2.1}.}
 \end{remark}

\begin{proof}
Clearly, $G(d)$ is continuous for all $d\ge0.$ In addition, $G(0)=0$ and
$G(d)\to\infty$ as $d\to\infty$ since (see \eqref{1.2})
$$G(d)\ge\int_{d/2}^d\int_{x-t}^{x+t}q(\xi)d\xi dt\ge
\frac{d}{2}\int_{x-d/2}^{x+d/2}q(\xi)d\xi\to\iy\quad\text{as}\quad d\to\iy.$$
Since we have, in addition,
\begin{equation}\label{2.7}
G'(d)=\int_{x-d}^{x+d}q(\xi)d\xi\ge0,\quad x\in R,\quad d\ge0,
\end{equation}
equation \eqref{2.3} has at least one solution $d_0>0.$

The obvious relations (see \eqref{2.7})
\begin{equation}\label{2.8}
1=G(d_0)=\int_0^{d_0}\int_{x-t}^{x+t}q(\xi)d\xi dt\le d_0\int_{x-d_0}^{x+d_0}
q(\xi)d\xi=d_0 G'(d_0)
\end{equation}
implies that $G'(d_0)>0$, and therefore $d_0$ is a unique root of \eqref{2.3}.
Denote it by $\hat d(x).$ {}From \eqref{2.8} and \lemref{lem2.2}, it follows
that $d(x)\le 2\hat d(x)$ since
$$2=2G(d(x))\le 2S(\hat d(x))\le S(2\hat d(x)).$$
The second inequality in \eqref{2.4} also follows from \lemref{lem2.2}:
$$2=2G(\hat d(x))\ge 2\int_{2\hat d(x)/3}^{\hat d(x)}\int_{x-t}^{x+t}q(\xi)d\xi
dt=\frac{2}{3}\hat d(x)\int_{x-2\hat d(x)/3}^{x+2\hat d(x)/3}
q(\xi)d\xi=S\left(\frac{2\hat d(x)}{3}\right).$$

Finally, the estimate \eqref{2.5} coincides with \eqref{2.8}, and it remains to
check \eqref{2.6}. Let us regard $\hat d(x)$ as an implicit function, i.e., as
the positive solution of the equation
\begin{equation}\label{2.9}
F(x,z)=\int_0^z\int_{x-t}^{x+t}q(\xi)d\xi dt-1=0.
\end{equation}
In  a neighborhood of the point $(x,\hat d(x))$, the function $F(x,z)$ is
continuous together with its partial derivatives
$$F_x'(x,z)=\int_0^z(q(x+t)-q(x-t))dt,\quad
F_z'(x,z)=\int_{x-z}^{x+z}q(\xi)d\xi.$$ In addition, according to \eqref{2.6},
we have
$$F_z'(x,z)\bigm|_{z=\hat d(x)}=\int_{x-\hat d(x)}^{x+\hat d(x)}q(\xi)d\xi\ge
\frac{1}{\hat d(x)}>0.$$ Hence $\hat d(x)$ is differentiable, and
\begin{equation}\label{2.10}
0=\hat d'(x)\int_{x-\hat d(x)}^{x+\hat d(x)}q(\xi)d\xi+\int_0^{\hat d(x)}
(q(x+t)-q(x-t))dt,\quad x\in R.
\end{equation}
{}From \eqref{2.10} and \eqref{2.5}, it now follows that
$$\frac{|\hat d'(x)|}{\hat d(x)}\le |\hat d'(x)|\int_{x-\hat d(x)}^{x+\hat
d(x)}q(\xi)d\xi=\left|\int_0^{\hat d(x)}(q(x+t)-q(x-t))dt\right|.$$
\end{proof}
\begin{cor}\label{cor2.5} If $q(x)\in H,$ then
\begin{equation}\label{2.11}
2k (x)|\hat d'(x)|\le 3c_2.
\end{equation}
\end{cor}

\begin{proof} {}From \eqref{1.12} and \eqref{2.4}, we get
\begin{equation}\label{2.12}
\hat d(x)\le\frac{3}{2} d(x)\le k (x)d(x),\quad x\in R.
\end{equation}
Therefore, according to \eqref{2.16}, \eqref{1.14} and \eqref{2.12}, we have
\begin{align*}
|\hat d'(x)|&\le\frac{3}{2}d(x)\left|\int_0^{\hat
d(x)}(q(x+t)-q(x-t))dt\right|\\
&\le\frac{3}{2}\frac{1}{k (x)}\left[k (x)d(x)\sup_{z\in[0,k (x)d(x)]}
\left|\int_0^z(q(x+t)-q(x-t))dt\right|\right]\le\frac{3c_2}{2k (x)}.
\end{align*}
\end{proof}

In what follows, we often us an obvious general assertion which, for
convenience, will be stated as a separate lemma.

\begin{lem}\label{lem2.6} Let $\varphi(x)$ and $\psi(x)$ be positive and
continuous for functions $f\in R.$ If there exists an interval $(a,b)$ such
that
\begin{equation}\label{2.13}
c^{-1}\varphi(x)\le\psi(x)\le c\varphi(x)\quad\text{for all}\quad x\notin
(a,b),
\end{equation}
then inequalities \eqref{2.13} remain true for all $x\in R$ (perhaps after
replacing with a larger constant).
\end{lem}

\begin{proof}
The function $f(x)=\frac{\psi(x)}{\varphi(x)}$ is continuous and positive for
$x\in [a,b].$ Hence its minimum $m$ and maximum $M$ on $[a,b]$ are finite
positive numbers. Let $\tilde c=\max\{c,m^{-1},M\}$ where $c$ is the constant
from \eqref{2.13}. Then $\tilde c^{-1}\varphi(x)\le\psi(x)\le\tilde
c\varphi(x)$ for $x\in R.$
\end{proof}

\begin{lem}\label{lem2.7}
Let $q(x)\in H$ and
\begin{equation}\label{2.14}
\omega(x)=[\omega^-(x),\omega^+(x)]=\left[x-\sqrt{k (x)}d(x),x+\sqrt{k
(x)}d(x)\right],\quad x\in R.
\end{equation}
Then there exists an absolute positive constant $c_3$ such that for all $x\in
R$ and $t\in\omega(x)$, the following inequalities hold:
\begin{equation}\label{2.15}
c_3^{-1}d(x)\le d(t)\le c_3d(x).
\end{equation}
\end{lem}

\begin{proof}
By \eqref{1.12}, there is $x_0\gg1$ such that $k (x)\ge36(c_1c_2)^2$ for
$|x|\ge x_0$ (see \eqref{1.13} and \eqref{1.14}).

In the following relations, we assume that $|x|\ge x_0,$\ $t\in\omega(x)$ and
use \eqref{2.11}, \eqref{1.13} and \eqref{2.4}:
\begin{align}
|\hat d(t)-\hat d(x)&=\left|\int_x^t\hat d'(\xi)d\xi\right|\le
\left|\int_x^t|\hat
d'(\xi)|d\xi\right|\le\frac{3c_2}{2}\left|\int_x^t\frac{d\xi}{k (\xi)}\right|\nonumber\\
&\le\frac{3c_1c_2}{2}\frac{|t-x|}{k (x)}\le\frac{3c_1c_2}{2}\frac{d(x)}{\sqrt{
k (x)}}\le\frac{3c_1c_2}{\sqrt{k (x)}}\hat d(x)\le\frac{\hat
d(x)}{2}\nonumber\\
&\Rightarrow\quad 2^{-1}\hat d(x)\le \hat d(t)\le 2\hat
d(x)\quad\text{for}\quad t\in\omega(x),\quad |x|\ge x_0.\label{2.16}
\end{align}

Denote
\begin{equation}\label{2.17}
\varphi(x)=\hat d(x),\qquad \psi_1(x)=\min_{t\in\omega(x)}\hat d(t),\qquad
\psi_2(x)=\max_{t\in\omega(x)}\hat d(t),\qquad x\in R.
\end{equation}
With the notation of \eqref{2.17}, inequalities \eqref{2.16} have the following
form:
\begin{equation}\label{2.18}
2^{-1}\varphi(x)\le\psi_1(x),  \psi_2(x)\le 2\varphi(x)\qquad\text{for}\quad
|x|\ge x_0.
\end{equation}
According to \eqref{2.18}, from \lemref{lem2.6} it follows that there exists a
constant $\tilde c$ such that
\begin{equation}\label{2.19}
\tilde c^{-1}\varphi(x)\le\psi_1(x),  \psi_2(x)\le\tilde
c\varphi(x)\qquad\text{for}\quad x\in R.
\end{equation}
The estimates \eqref{2.19} immediately imply the inequalities
\begin{equation}\label{2.20}
\tilde c^{-1}\hat d(x)\le\hat d(t)\le\tilde c\hat d(x)\qquad\text{for}\quad
t\in\omega(x),\quad x\in R.
\end{equation}
The relations \eqref{2.15} with $c_3=3\tilde c$ follow from \eqref{2.20} and
\eqref{2.4}:
$$\frac{d(t)}{d(x)}=\frac{d(t)}{\hat d(t)}\cdot \frac{\hat d(t)}{\hat d(x)}\cdot
\  \frac{\hat d(x)}{d(x)}\le 2\cdot\tilde c\cdot \frac{3}{2}=3\tilde
c\qquad\text{for}\quad t\in\omega(x),\quad x\in R,$$
$$\frac{d(t)}{d(x)}=\frac{d(t)}{\hat d(t)}\cdot \frac{\hat d(t)}{\hat d(x)}\cdot
 \frac{\hat d(x)}{d(x)}\ge
\frac{2}{3}\cdot\frac{1}{\tilde c}\cdot \frac{1}{2}=\frac{1}{3\tilde
c}\qquad\text{for}\quad t\in\omega(x),\quad x\in R.$$
\end{proof}
\begin{lem}\label{lem2.8}
Under condition \eqref{1.2}, we have
\begin{equation}\label{2.21}
\lim_{x\to-\infty}(x+d(x))=-\infty,\qquad\lim_{x\to\infty}(x-d(x))=\infty.
\end{equation}
\end{lem}

\begin{proof} The equalities in \eqref{2.21} are checked in a similar way. Let
us prove, for example, the second one. We show that
\begin{equation}\label{2.22}\varliminf_{x\to\infty}(x-d(x))=\iy.
\end{equation}
Assume the contrary. Then there exists a number $a\in R$ and a sequence
$\{x_n\}_{n=1}^\iy$ such that
\begin{equation}\label{2.23}x_n-d(x_n)\le a\quad\text{for}\quad n\in N\qquad\text{and}\qquad
x_n\to\iy\quad\text{as}\quad n\to\iy.
\end{equation}
 {}From \eqref{2.23} it follows that
there is $n_0\gg1$ such that for all $n\ge n_0$, the following inequalities
hold:
\begin{equation}\label{2.24}
d(x_n)\ge x_n-a=x_n\left(1-\frac{a}{x_n}\right)\ge\frac{x_n}{2},\quad n\ge n_0.
\end{equation}
Then using \eqref{2.23}, \eqref{2.24} and \eqref{1.11}, we get
\begin{equation}\label{2.25}2=d(x_n)\int_{x_n-d(x_n)}^{x_n+d(x_n)}q(t)dt\ge\frac{x_n}{2}\int_a^{x_n}q(t)dt\quad
\Rightarrow\quad\frac{4}{x_n}\ge\int_a^{x_n}q(\xi)d\xi.
\end{equation}

Clearly, \eqref{2.23} and \eqref{2.25} contradict \eqref{1.2}, which leads to
\eqref{2.22}. But this implies the statement of the lemma because
$$\iy=\varliminf_{x\to\iy}(x-d(x))\le\varlimsup_{x\to\iy}(x-d(x))\le\iy\quad\Rightarrow
\quad \varliminf_{x\to\iy}(x-d(x))=\varlimsup_{x\to\iy}(x-d(x))=\iy.$$
\end{proof}

\section{Main asymptotic formula}

In this section, our goal is to prove the following assertion.

\begin{lem}\label{lem3.1}
Suppose that condition \eqref{1.2} holds and
\begin{equation}\label{3.1}
\lim_{|x|\to\iy}\rho'(x)=0.
\end{equation}
Then we have
\begin{equation}\label{3.2}
\rho(x)=\frac{d(x)}{2}(1+\varepsilon(x)),\qquad\lim_{|x|\to\iy}\varepsilon(x)=0,
\end{equation}
and the following relations hold:
\begin{equation}\label{3.3}
|\varepsilon(x)|\le ch(x),\qquad h(x)\to0\qquad\text{as}\quad |x|\to\iy.
\end{equation}
\end{lem}

Here
\begin{equation}\label{3.4}
h(x)=\sup_{t\in\Delta(x)}|\rho'(t)|,\qquad
\Delta(x)=[\Delta^-(x),\Delta^+(x)]=[x-d(x),x+d(x)],\qquad x\in R.
\end{equation}
Denote by $z_0$ a point on the number axis such that (see \eqref{3.1})
\begin{equation}\label{3.5}|\rho'(x)|\le 10^{-3}\qquad\text{for all}\qquad |x|\ge z_0.
\end{equation}
Suppose that in addition to \eqref{1.2}, \eqref{3.1}, we have $d_0<\iy$ (see
\eqref{1.23}). Then equality \eqref{3.2} holds, and
\begin{equation}\label{3.6}|\varepsilon(x)|\le 18h(x)\qquad\text{for}\qquad
|x|\ge z_1,\qquad z_1\doe z_0+d_0+1.
\end{equation}

To prove \lemref{lem3.1}, we need the following auxiliary assertions.

\begin{lem}\label{lem3.2}
For $x\in R$ the following relations hold:
\begin{equation}\label{3.7}
|\rho'(x)|<1,
\end{equation}
\begin{equation}\label{3.8}
\frac{v'(x)}{v(x)}=\frac{1+\rho'(x)}{2\rho(x)},\qquad
\frac{u'(x)}{u(x)}=-\frac{1-\rho'(x)}{2\rho(x)}.
\end{equation}
\end{lem}

\begin{proof}
Let us show that (see \eqref{1.3})
\begin{equation}\label{3.9}
v'(x)>0,\qquad u'(x)<0\qquad\text{for}\qquad x\in R.
\end{equation}

For a given $x\in R,$ by \eqref{1.2} there exists $a\in (-\iy,x]$ such that
$$\int_a^xq(t)dt>0.$$
Then from \eqref{1.2} and \eqref{1.3} it follows that
$$v'(x)=v'(a)+\int_z^xq(t)v(t)dt\ge\int_a^xq(t)v(t)dt\ge
v(a)\int_a^xq(t)dt>0.$$

The second inequality from \eqref{3.9} can be checked in a similar way. To
prove \eqref{3.8}, it suffices to differentiate \eqref{1.9}. Inequality
\eqref{3.7} follows from \eqref{3.8} and \eqref{3.9}.
\end{proof}

\begin{lem}\label{lem3.3}
For $x\in R$, we have \begin{equation}\label{3.10}
\frac{1+\rho'(\Delta^+(x))}{1-\rho'(\Delta^+(x))}\cdot\frac{1-\rho'(\Delta^-(x))}{
1+\rho'(\Delta^-(x))}=\exp\left(4\int_{\Delta(x)}\frac{q(t)\rho(t)dt}{1-\rho'(t)^2}
-\int _{\Delta(x)}\frac{dt}{\rho(t)}\right).
\end{equation}
\end{lem}

Here $\Delta(x)=[\Delta^-(x),\Delta^+(x)]=[x-d(x),x+d(x)].$

\begin{proof}
{}From \eqref{3.9} and \eqref{1.1} for $t\in R$, it follows that
$$v''(t)=q(t)v(t)\quad\Rightarrow\quad
\frac{v''(t)}{v'(t)}=q(t)\frac{v(t)}{v'(t)}\quad\Rightarrow\quad
ln\frac{v'(\Delta^+(x))}{v'(\Delta^-(x))}=\int_{\Delta(x)}\frac{q(t)v(t)dt}{v'(t)},$$
$$u''(t)=q(t)u(t)\quad\Rightarrow\quad
\frac{u''(t)}{u'(t)}=q(t)\frac{u(t)}{u'(t)}\quad\Rightarrow\quad
ln\frac{u'(\Delta^+(x))}{u'(\Delta^-(x))}=\int_{\Delta(x)}\frac{q(t)u(t)dt}{u'(t)}.$$
These inequalities imply
\begin{equation}\label{3.11}
\frac{v'(\Delta^+(x))}{v'(\Delta^-(x))}\cdot\frac{u'(\Delta^-(x))}{u'(\Delta^+(x))}=
\exp\left(\int_{\Delta(x)}q(t)\left(\frac{v(t)}{v'(t)}-\frac{u(t)}{u'(t)}\right)dt\right),
\quad x\in R.
\end{equation}

When substituting \eqref{3.8} into \eqref{3.11}, we get
\begin{equation}\label{3.12}
\frac{1+\rho'(\Delta^+(x))}{1-\rho'(\Delta^+(x))}\cdot
\frac{1-\rho'(\Delta^-(x))}{1+\rho'(\Delta^-(x))}\cdot\frac{v(\Delta^+(x))}{u(\Delta^+
(x))}\cdot
\frac{u(\Delta^-(x))}{v(\Delta^-(x))}=\exp\left(4\int_{\Delta(x)}\frac{q(t)\rho(t)dt}{
1-\rho'(t)^2}\right).
\end{equation}

Furthermore, according to \eqref{1.9}  we have
\begin{equation}\label{3.13}
\frac{v(\Delta^+(x))}{u(\Delta^+(x))}=\exp\left(\int_{x_0}^{\Delta^+(x)}\frac{dt}{\rho(t)}
\right),\qquad
\frac{v(\Delta^-(x))}{u(\Delta^-(x))}=\exp\left(\int_{x_0}^{\Delta^-(x)}\frac{dt}{\rho
(t)}\right),\qquad x\in R.
\end{equation}
To prove \eqref{3.10}, it remains to substitute \eqref{3.13} into \eqref{3.12}.
\end{proof}

\begin{lem}\label{lem3.4}
Suppose that condition \eqref{1.2} holds. Then
\begin{equation}\label{3.14}
\rho(t)\le\frac{5}{2}d(x)\qquad\text{for}\qquad
t\in\Delta(x)=[x-d(x),x+d(x)],\qquad x\in R.
\end{equation}
\end{lem}

\begin{proof} By Lagrange's formula,
\begin{equation}\label{3.15}
\rho(t)=\rho(x)+\rho'(\xi)(t-x),\qquad t\in\Delta(x),\qquad x\in R.
\end{equation}
The point $\xi$ in \eqref{3.15} lies between $t$ and $x$. Then \eqref{3.15},
together with \eqref{3.7} and \eqref{3.10}, lead to \eqref{3.14}:
$$\rho(t)\le\rho(x)+|\rho'(\xi)|\
|t-x|\le\rho(x)+d(x)\le\frac{3}{2}d(x)+d(x)=\frac{5}{2}d(x).$$
\end{proof}

In the sequel, we assume that conditions \eqref{1.2} and \eqref{3.1} hold and
do not mention them in the statements.

\begin{lem}\label{lem3.5}
For all $|x|\gg1,$ the following inequalities hold:
\begin{equation}\label{3.16}
\left|4\int_{\Delta(x)}\frac{q(t)\rho(t)dt}{1-\rho'(t)^2}-8\frac{\rho(x)}{d(x)}\right|
\le 8,0201h(x).
\end{equation}
In addition, if $d_0<\iy$ (see \eqref{1.23}), then \eqref{3.16} holds for all
$|x|\ge z_1$ (see \eqref{3.6}).
\end{lem}

\begin{proof} In the following transformations, we use the definition of $d(x)$
(see \eqref{1.11}):
\begin{align}
4\int_{\Delta(x)}&\frac{q(t)\rho(t)dt}{1-\rho'(t)^2}
=4\int_{\Delta(x)}q(t)\rho(t)dt+4\int
_{\Delta(x)}\frac{q(t)\rho(t)\rho'(t)^2}{1-\rho'(t)^2}dt\nonumber\\
&=4\rho(x)\int_{\Delta(x)}q(t)dt+4\int_{\Delta(x)}q(t)(\rho(t)-\rho(x))dt+4\int_{\Delta
(x)}\frac{q(t)\rho(t)\rho'(t)^2}{1-\rho'(t)^2}dt\nonumber\\
&=\frac{8\rho(x)}{d(x)}+4\int_{\Delta(x)}q(t)(\rho(t)-\rho(x))dt+4\int_{\Delta(x)}
\frac{q(t)\rho(t)\rho'(t)^2}{1-\rho'(t)^2}dt,\ x\in R.\label{3.17}
\end{align}

Below, in the estimate of the first integral of \eqref{3.17}, we use
\eqref{3.15} and the definitions of $h(x)$,\ $\Delta(x)$ and $d(x)$ (see
\eqref{3.4}, \eqref{1.11}):
\begin{align}
4\left|\int_{\Delta(x)}q(t)(\rho(t)-\rho(x))dt\right|&\le
4\int_{\Delta(x)}q(t)|\rho(t)-\rho(x)|dt=4\int_{\Delta(x)}q(t)|\rho'(\xi)|\
|t-x|dt\nonumber\\
&\le 4h(x)d(x)\int_{\Delta(x)}q(t)dt=8h(x),\qquad x\in R.\label{3.18}
\end{align}

Let us estimate the second integral from \eqref{3.17}. {}From \eqref{2.21}, it
follows that there is $\tilde z_0\gg z_0$ (see \eqref{3.5}) such that
\begin{equation}\label{3.19}
\Delta(x)\cap[-z_0,z_0]=\emptyset\qquad\text{for}\quad |x|\gg\tilde z_0.
\end{equation}
In particular, if $d_0<\infty$ (see \eqref{1.23}), then one can set $\tilde
z_0:=z_1=z_0+d_0+1$ (see \eqref{3.6}). Indeed, with such a choice of $\tilde
z_0$, we have
\begin{enumerate}\item[1)] if $x\le-\tilde z_0\ \Rightarrow\ x+d(x)\le
-\tilde z_0+d(x)=-z_0-1+d(x)-d_0<-z_0\ \Rightarrow\ \eqref{3.19}$ \item[2)] if
$x\ge-\tilde z_0\ \Rightarrow\ x-d(x)\ge\tilde z_0-d(x)=z_0+1+d_0-d(x)>z_0\
\Rightarrow\ \eqref{3.19}.$
\end{enumerate}

Below, for $|x|\ge\tilde z_0,$ we use \eqref{3.19}, \eqref{3.14}, \eqref{3.5}
and \eqref{1.11}:
\begin{align}
0&\le 4\int_{\Delta(x)}\frac{q(t)\rho(t)\rho'(t)^2dt}{1-\rho'(t)^2}\le
\frac{4\cdot 10^{-3}}{1-10^{-6}}h(x)\int_{\Delta(x)}q(t)\rho(t)dt
\nonumber\\&\le
\frac{4h(x)}{10^3-10^{-3}}\frac{5}{2}d(x)\int_{\Delta(x)}q(t)dt=\frac{20h(x)}
{10^3-10^{-3}}\le 0.0201h(x).\label{3.20}
\end{align}
{}From \eqref{3.20} and \eqref{3.18}, we get \eqref{3.16}
\end{proof}

\begin{lem}\label{lem3.6}
For all $|x|\gg1,$ we have
\begin{equation}\label{3.21}
\left|\int_{\Delta(x)}\frac{dt}{\rho(t)}-2\frac{d(x)}{\rho(x)}\right|\le
32.16h(x).\end{equation} In addition, if $d_0<\infty$ (see \eqref{1.23}), then
\eqref{3.21} holds for all $|x|\ge z_1$ (see \eqref{3.6}).
\end{lem}

\begin{proof}
Let $\tilde z_0$ be the number from \lemref{lem3.5}. In the following
transformation, we use \eqref{3.15}:
\begin{equation}\label{3.22}
\int_{\Delta(x)}\frac{dt}{\rho(t)}=\frac{1}{\rho(x)}\int_{\Delta(x)}
\frac{\rho(x)dt}{\rho(x
)+\rho'(\xi)(t-x)}=\frac{1}{\rho(x)}\int_{\Delta(x)}\frac{dt}{1+\rho'(\xi)\frac{t-x}{
\rho(x)}}.
\end{equation}

Consider the integrand in \eqref{3.22}. Let us check the estimate
\begin{equation}\label{3.23}
|\gamma(x,\xi,t)|\le 4\cdot 10^{-3},\quad \gamma(x,\xi,t)\doe
\rho'(\xi)\frac{t-x}{\rho(x)},\quad \xi,t\in\Delta(x),\quad |x|\ge \tilde
z_0.\end{equation} Indeed, for $|x|\ge\tilde z_0$ from \eqref{1.10}, it follows
that
$$|\gamma(x,\xi,t)|=|\rho'(\xi)|\frac{|t-x|}{\rho(x)}\le
10^{-3}\frac{dx}{\rho(x)}\le 4\cdot 10^{-3}\ \Rightarrow\ \eqref{3.23}.$$

Below, for $|x|\ge\tilde z_0$, we use \eqref{3.23}, \eqref{1.10} and the
definition of $h(x)$ (see \eqref{3.4}):
\begin{align}
\frac{|\gamma(x,\xi,t)|}{|1+\gamma(x,\xi,t)|}&\le
\frac{|\gamma(x,\xi,t)|}{1-|\gamma(x,\xi,t)|}\le\frac{h(x)}{1-4.10^{-3}}\frac{d(x)}{\rho(x)}
\nonumber\\ &\le \frac{4000}{996}h(x)\le 4.02h(x),\quad |x|\ge\tilde
z_0.\label{3.24}\end{align}

To finish the proof of \eqref{3.21}, it remains to apply \eqref{3.22},
\eqref{3.24} and \eqref{1.10} for $|x|\ge\tilde z_0:$
\begin{align*}
 \left|\int_{\Delta(x)}\frac{dt}{\rho(t)}-2\frac{d(x)}{\rho(x)}\right|&=\left|\frac{1}
{\rho(x)}\int_{\Delta(x)}\frac{dt}{1+\gamma(x,\xi,t)}-2\frac{d(x)}{\rho(x)}\right|
\\ &=\frac{1}{\rho(x)}\left|\int_{\Delta(x)}\left(\frac{1}{1+\gamma(x,\xi,t)}-1\right)dt
\right|
 \\
 &
\le
\frac{1}{\rho(x)}\int_{\Delta(x)}\frac{|\gamma(x,\xi,t)|dt}{|1+\gamma(x,\xi,t)|}\le
8.04 h(x)\frac{d(x)}{\rho(x)}\le 32.16h(x).\end{align*}
\end{proof}

\begin{proof}[Proof of \lemref{lem3.1}]
Throughout the sequel, we assume $|x|\ge\tilde z_0$ where $\tilde z_0$ is the
number from \lemref{lem3.5}. Consider \eqref{3.10}. In the following estimates,
we use \eqref{3.5} and the definition of $h(x)$ (see \eqref{3.4}):
\begin{gather}
 \left|\frac{1+\rho'(\Delta^+(x))}{1-\rho'(\Delta^+(x))}\cdot\frac{1-\rho'(\Delta^-(x))}
{1+\rho'(\Delta^-(x))}-1\right|
=\frac{2|\rho'(\Delta^+(x))-\rho'(\Delta^-(x))|}{
|+\Delta'(\Delta^-(x))-\rho'(\Delta^+(x))-\rho'(\Delta^-(x))\rho'(\Delta^+(x))|}\nonumber \\
   \qquad\quad\qquad\qquad\le\frac{4h(x)}{1-2.10^{-3}-10^{-6}}\le 4.009h(x) \nonumber\\
 \Rightarrow\ \frac{1+\rho'(\Delta^+(x))}{1-\rho'(\Delta^-(x))}\
\frac{1-\rho'(\Delta^-(x))}{1+\rho '(\Delta^-(x))}=1+\delta_1(x),\quad
|\delta_1(x)|\le 4.009h(x),\quad |x|\ge\tilde z_0. \label{3.25}
\end{gather}

Below, in the transformation of the exponent in \eqref{3.10}, we use
inequalities \eqref{3.16} and
 \eqref{3.21}:
 \begin{gather}
  \left|\left(4\int_{\Delta(x)}\frac{q(t)\rho(t)dt}{1-\rho'(t)^2}-\int_{
 \Delta(x)}\frac{dt}{\rho(t)}\right)-\left(\frac{8\rho(x)}{d(x)}-\frac{2d(x)}{\rho(x)}
 \right)\right|\nonumber\\
 \qquad \le\left| 4\int_{\Delta(x)}\frac{q(t)\rho(t)dt}{1-\rho'(t)^2}-
\frac{8\rho(x)}{d(x)}\right|+\left| \int_{
 \Delta(x)}\frac{dt}{\rho(t)}
-\frac{2d(x)}{\rho(x)}\right|\le(8,0201+32.16)h(x)\le 40.2h(x)\nonumber\\
 \Rightarrow\
4\int_{\Delta(x)}\frac{q(t)\rho(t)dt}{1-\rho'(t)^2}-\int_{\Delta(x)}\frac{dt}{\rho(t)}
= \frac{8\rho(x)}{d(x)}-\frac{2d(x)}{\rho(x)}+\delta_2(x),\nonumber\\
 |\delta_2(x)|\le 40.2h(x),\quad |x|\ge\tilde
z_0.\label{3.26}
 \end{gather}
Thus, (see \eqref{3.25} and \eqref{3.26}) equality \eqref{3.10} is reduced to
\begin{equation}\label{3.27}
1+\delta_1(x)=\exp\left(8\frac{\rho(x)}{d(x)}-2\frac{d(x)}{\rho(x)}+\delta_2(x)\right),
\qquad |x|\ge\tilde z_0.
\end{equation}
{}From Lagrange's formula, \eqref{3.5} and \eqref{3.25} it follows that
\begin{gather}
ln(1+\delta_1(x))=\frac{\delta_1(x)}{1+\xi},\qquad
\xi\in(-|\delta_1(x)|,|\delta(x)|)\nonumber\\ \Rightarrow\
ln(1+\delta_1(x))=\delta_3(x),\quad
|\delta_3(x)|\le\frac{|\delta_1(x)|}{1-|\delta_1(x)|}\le\frac{4.009h(x)}{1-4.009\cdot
10^{-3}}<4.03h(x).\label{3.28}
\end{gather}
According to \eqref{3.27} and \eqref{3.28}, we now obtain
\begin{gather}
\delta_3(x)=\frac{8\rho(x)}{d(x)}-\frac{2d(x)}{\rho(x)}+\delta_2(x),\quad
|x|\ge\tilde z_0\nonumber\\
\Rightarrow\ \frac{8\rho(x)}{d(x)}-\frac{2d(x)}{\rho(x)}=\delta_4(x),\quad
|\delta_4(x)|\le|\delta_2(x)|+|\delta_3(x)|\le 44.23h(x),\quad|x|\ge\tilde
z_0.\label{3.29}
\end{gather}

Let us rewrite \eqref{3.29} in the following way:
\begin{equation}\label{3.30}
\rho(x)^2=\frac{d(x)^2}{4}\left(1+\delta_4(x)\frac{\rho(x)}{2d(x)}\right),\quad
|\delta_4(x)|\le 44.23h(x),\quad|h|\ge\tilde z_0.
\end{equation}
Denote
\begin{equation}\label{3.31}
\alpha(x)=\delta_4(x)\frac{\rho(x)}{2d(x)},\qquad |x|\ge \tilde z_0.
\end{equation}

Below, in the estimate of $|\alpha(x)|,$ we use \eqref{3.30} and \eqref{1.10}:
\begin{equation}\label{3.32}
|\alpha(x)|\le|\delta_4(x)|\frac{\rho(x)}{2d(x)}\le
44.23\cdot\frac{3}{4}h(x)\le 33.2h(x)\le 0.0332.\end{equation} Therefore, from
\eqref{3.30}, \eqref{3.31} and \eqref{3.32}, we get
\begin{equation}\label{3.33}
\rho(x)=\frac{d(x)}{2}\sqrt{1+\alpha(x)},\qquad |x|\ge\tilde z_0.
\end{equation}
Furthermore, since
\begin{equation}\label{3.34}
\sqrt{1+\nu}=1+\frac{\nu}{2}-\frac{1}{2}\left(\frac{\nu}{1+\sqrt{1+\nu}}\right)^2\qquad\text{for}
\quad 1+\nu\ge0,
\end{equation}
from \eqref{3.32} and \eqref{3.34}, we get
\begin{equation}\label{3.35}
\sqrt{1+\alpha(x)}=1+\varepsilon(x),\quad\
|\varepsilon(x)|\le\frac{|\alpha(x)|}{2}+\frac{1}{2}\left(\frac{\alpha(x)}{1+\sqrt{1+\alpha(
x)}}\right)^2,\quad\ |x|\ge\tilde z_0.
\end{equation}

In the following estimate of $|\varepsilon(x)|$, we use \eqref{3.35},
\eqref{3.32}, \eqref{3.30} and \eqref{1.10}:
\begin{align}
|\varepsilon(x)|&\le\frac{|\alpha(x)|}{2}+\frac{|\alpha(x)|^2}{2}=|\alpha(x)|\frac{
|1+|\alpha(x)|}{2}\le
\frac{1+0.0332}{2}|\alpha(x)|=0.5166|\alpha(x)|\nonumber\\
&=0.5166|\delta_4(x)|\frac{\rho(x)}{2d(x)}\le 44.23\cdot \frac{3}{4}\cdot
0.5166h(x)<18h(x),\qquad |x|\ge z_0.\label{3.36} \end{align} \lemref{lem3.1}
now follows from \eqref{3.33}, \eqref{3.35} and \eqref{3.36}.
\end{proof}

\section{Proof of the main result}

In this section we finish the proof of formula \eqref{1.19}. Note that this
part more or less coincides with the corresponding fragment of \cite{6} and is
reproduced here, with minor changes, only for the sake of completeness.

\begin{lem}\label{lem4.1} For $x\in R,$ we have the inequality
\begin{equation}\label{4.1}
|\rho'(x)|\le|\varkappa(x)-1|,\qquad
\varkappa(x)\doe\frac{v'(x)}{v(x)}\cdot\frac{u(x)}{|u'(x)|}.
\end{equation}
\end{lem}

\begin{proof}
{}From \eqref{1.3}, \eqref{1.4} and \eqref{3.9}, it follows that
$$|\rho'(x)|=|v'(x)u(x)+v(x)u'(x)|=|u'(x)|v(x)|\varkappa(x)-1|<|\varkappa(x)-1|.$$
\end{proof}

\begin{lem}\label{lem4.2}
For $x\in R,$ the formula
\begin{equation}\label{4.2}
y(t)=v'(x)u(t)-u'(x)v(t),\qquad t\in R\end{equation} determines the solution of
the Cauchy problem
\begin{equation}\label{4.3}
y''(t)=q(t)y(t),\qquad t\in R,
\end{equation}
\begin{equation}\label{4.4}
y(t)\bigm|_{t=x}=1,\qquad y'(t)\bigm|_{t=x}=0. \end{equation} In addition, the
following inequalities hold:
\begin{equation}\label{4.5}
y'(t)\le 0\qquad\text{for}\quad t\le x;\qquad y'(t)\ge0\qquad\text{for}\quad
t\ge x .\end{equation}
\end{lem}

\begin{proof}
Let us check \eqref{4.5} for $t\ge x.$ (The other assertions of the lemma
immediately follows from the properties of the PFSS $\{u(x),v(x)\}$ of equation
\eqref{1.1} (see \eqref{1.4}.) Let us show that $y(t)>0$ for $t>x.$ If this is
not the case, let $x_0$ be the smallest positive root of the equation $y(t)=0$
\ $(x_0>0$ because of \eqref{4.4}). Then $y'(x_0)\le 0.$ Indeed, if
$y'(x_0)>0,$ then $y(t)<0$ for $t<x_0$ because $y(x_0)=0.$ But then \eqref{4.4}
implies that the equation $y(t)=0$ has a root in the interval $(0,x_0)$ which
contradicts the definition of $x_0.$ Thus $y'(x_0)\le 0.$ On the other hand,
from \eqref{4.3} it follows that
$$y'(x_0)=\int_0^{x_0}q(\xi) y(\xi)d\xi\ge0\quad\Rightarrow\quad y'(x_0)=0.$$
Hence $y(t)=0$ because $y(x_0)=y'(x_0)=0.$ Contradiction.

Since $y(t)>0$ for $\ge x,$ according to \eqref{1.2} and
\eqref{4.3}--\eqref{4.4} we get
$$y'(t)=\int_0^tq(\xi)y(\xi)d\xi\ge0\qquad\text{for}\quad t\ge x.$$
The case $t\le x$ is treated in a similar way.
\end{proof}

Let $q(x)\in H.$ Let us introduce the functions (see \eqref{2.14})
\begin{equation}\label{4.6}
\tilde u(t)=y(t)\int_t^{\omega^+(x)}\frac{d\xi}{y(\xi)^2},\quad \tilde
v(t)=y(t)\int_{\omega^-(x)}^t\frac{d\xi}{y(\xi)^2},\quad t\in\omega(x),\quad
x\in R.
\end{equation}
In \eqref{4.6}, we assume that $y(t)$ is the solution of the problem
\eqref{4.3}--\eqref{4.4}.

\begin{lem}\label{lem4.3}
The functions \eqref{4.6} are solutions of equation \eqref{4.3} and satisfy the
relations
\begin{equation}\label{4.7}
\tilde u(\omega^+(x))=\tilde v(\omega^-(x))=0,\quad \tilde u(t)\ge0,\quad
\tilde v(t)\ge0\qquad\text{for}\quad t\in\omega(x),
\end{equation}
\begin{equation}\label{4.8}
\tilde v'(t)\tilde u(t)-  \tilde u'(t)  \tilde
v(t)=\int_{\omega(x)}\frac{d\xi}{y(\xi)^2},\quad t\in\omega(x).
\end{equation}
\end{lem}

\begin{proof}
The relations \eqref{4.7} are obvious. Equality \eqref{4.8} is checked by a
straightforward calculation.
\end{proof}
\begin{lem}\label{lem4.4}
For $x\in R,$ we have the equalities
\begin{equation}\label{4.9}
\frac{v(x)}{v'(x)}=\tilde v(x)\left[1-\frac{v(\omega^-(x))}{v(x)}\
\frac{1}{y(\omega^-(x))}\right]^{-1},
\end{equation}
\begin{equation}\label{4.10}
\frac{u(x)}{|u'(x)|}=\tilde u(x)\left[1-\frac{u(\omega^+(x))}{u(x)}\
\frac{1}{y(\omega^+(x))}\right]^{-1}.
\end{equation}
\end{lem}

Here $y(\cdot)$ is the solution of problem \eqref{4.3}--\eqref{4.4}.

\begin{proof}
The equalities \eqref{4.9}--\eqref{4.10} follow from \eqref{1.9}, \eqref{4.2}
and \eqref{4.6}. For example,
\begin{align*}
\tilde v(x)&=\int_{\omega^-(x)}^x\frac{d\xi}{[v'(x)u(\xi)-u'(x)v(\xi)]^2}=
\int_{\omega^-(x)}^x \frac{1}{\rho(\xi)}\
\frac{\exp\left(-\int_{x_0}^\xi\frac{ds}{\rho(s)}\right)d\xi}{\left[v'(x)\exp\left(
-\int_{x_0}^\xi\frac{ds}{\rho(s)}\right)-u'(x)\right]^2}\\
 &=\int_{\omega^-}^x\frac{1}{v'(x)}d\left[v'(x)\exp\left(-\int_{x_0}^\xi
 \frac{ds}{\rho(s)}
 \right)-u'(x)\right]^{-1}\\ &=\frac{1}{v'(x)}\left[v'(x)\exp\left(-\int_{x_0}^\xi\frac{ds}
{\rho(s)}\right)-u'(x)\right]^{-1}\bigm|_{\omega^-(x)}^x\\
 &=\frac{1}{v'(x)}\
 \frac{v(\xi)}{v'(x)u(\xi)-u'(x)v(\xi)}\bigm|_{\omega^-(x)}^x=\frac{1}{v'(x)}\
 \frac{v(\xi)}{y(\xi)}\bigm|_{\omega^-(x)}^x=\frac{1}{v'(x)}\left[v(x)-\frac{
 v(\omega^-(x))}{y(\omega^-(x))}\right]\\
 &=\frac{v(x)}{v'(x)}\left[1-\frac{v(\omega^-(x))}{v(x)}\
\frac{1}{y(\omega^-(x))}\right]\ \Rightarrow\ \eqref{4.9}.
\end{align*}

The equality \eqref{4.10} is checked in a similar way.
\end{proof}

\begin{lem}\label{lem4.5}
Suppose $q(x)\in H.$ Then for $x\in R$ we have the equality (see \eqref{2.15}:
\begin{equation}\label{4.11}
\varkappa(x)=\frac{v'(x)}{v(x)}\ \frac{u(x)}{|u'(x)|}=\frac{\tilde u(x)}{\tilde
v(x)}(1+\nu(x)),\quad |\nu(x)|\le\sqrt{6c_3}\exp\left(-\frac{\sqrt{k
(x)}}{3c_3}\right).
\end{equation}
\end{lem}

\begin{proof}
By \lemref{lem4.4}, to prove \eqref{4.11} it is enough to show that for $x\in
R,$ the following inequality holds:
\begin{equation}\label{4.12}
\max\left\{\frac{v(\omega^-(x))}{v(x)}\ \frac{1}{y(\omega^-(x))},\quad
\frac{u(\omega^+(x))}{u(x)}\
\frac{1}{y(\omega^+(x))}\right\}\le\sqrt{6c_3}\exp\left(-\frac{\sqrt{k
(x)}}{3c_3}\right).
\end{equation}
By \eqref{4.5}, \eqref{1.9}, \eqref{1.10} and \eqref{2.15}, we have
\begin{align*}
\frac{u(\omega^+(x))}{u(x)}\ \frac{1}{y(\omega^+(x))}&\le
\frac{\omega^+(x))}{u(x)}\\
&=\sqrt{\frac{\rho(\omega^+(x))}{d(\omega^+(x))}\cdot
\frac{d(\omega^+(x))}{d(x)}\cdot\frac{d(x)}{\rho(x)}}\exp\left(-\frac{1}{2}\int_x^{\omega^+(x)}
\frac{d(\xi)}{\rho(\xi)}\cdot \frac{d(x)}{d(\xi)}\cdot
\frac{d\xi}{d(x)}\right)\\
 &\le \sqrt{\frac{3}{2}\cdot c_3\cdot
4}\exp\left(-\frac{1}{2}\int_x^{\omega_+(x)}\frac{2}{3}\ \frac{1}{c_3}\
 \frac{d\xi}{d(x)}\right)=\sqrt{6c_3}\exp\left(-\frac{\sqrt{k (x)}}{3c_3}\right).
\end{align*}

The second inequality of \eqref{4.12} is checked in a similar way.
\end{proof}

To study $\tilde u(x)/\tilde v(x)$, let us look  at the solution $y(t)$ of
problem \eqref{4.3}--\eqref{4.4} more closely than in \lemref{lem4.2}. Suppose
$q(x)\in H.$ Denote
\begin{equation}\label{4.13}
\chi(x)=[0,\sqrt{k (x)}d(x)],\qquad x\in R.
\end{equation}

In problem \eqref{4.3}--\eqref{4.4}, we change variables:
\begin{equation}\label{4.14}
y_1(x)=y(x-z),\qquad z\in \chi(x),
\end{equation}
\begin{equation}\label{4.15}
y_2(x)=y(x+z),\qquad z\in \chi(x).
\end{equation}
It is easy to see that $y_1(z)$ and $y_2(z)$ are solutions of the following
Cauchy problems, respectively:
\begin{equation}\label{4.16}
y_1''=q(x-z)y_1(x),\qquad y_1(0)=1,\qquad y_1'(0)=0,
\end{equation}
\begin{equation}\label{4.17}
y_2''=q(x+z)y_1(x),\qquad y_2(0)=1,\qquad y_2'(0)=0.
\end{equation}

\begin{lem}\label{lem4.6}
Suppose $q(x)\in H,$ and let $t_0$ be a positive number such that $k (x)\ge
64c_2^2$ for all $|x|\ge t_0$ (see \eqref{1.12}, \eqref{1.14}). Then for
$|x|\ge t_0$ and $z\in\chi(x)$ (see \eqref{4.13}), the following relations
hold:
\begin{equation}\label{4.18}
\frac{y_2(z)}{y_1(z)}=\frac{y(x+z)}{y(x-z)}=1+\gamma(z),\qquad |\gamma(z)|\le
1.2F(x)\le\frac{1.2c_2}{\sqrt{k (x)}}.
\end{equation}
\end{lem}

\begin{proof} Let us introduce some notation:
\begin{equation}\label{4.19}
\beta(z)=\frac{y_2(z)}{y_1(z)},\quad
\varphi(z)=\int_0^z(q(x+\xi)-q(x-\xi))d\xi,\quad\psi(z)=\max_{t\in[0,z]}|\varphi(t)|,\quad
z\in\chi(x).
\end{equation}
By \eqref{4.5}, we get $y_1'(z)\ge0,$\ $y_2'(z)\ge0$ for $z\ge0.$ Therefore, we
also have $[y_1(z)\cdot y_2(z)]'\ge0$ for $z\ge0.$ Integrating by parts, we get
\begin{align}
|\beta'(z)|&=\left|\frac{d}{dz}\left(\frac{y_2(z)}{y_1(z)}\right)\right|=
\frac{1}{y_1(z)^2}\left|\int_0^z(q(x+\xi)
-q(x-\xi))y_1(\xi)y_2(\xi)d\xi\right|\nonumber\\
&\le
\frac{y_2(z)}{y_1(z)}|\varphi(z)|+\frac{1}{y_1(z)^2}\left|\int_0^z\varphi(\xi)[y_1(\xi)y_2
(\xi)]'d\xi\right| \nonumber\\
 &\le
\frac{y_2(z)}{y_1(z)}\psi(z)+\frac{\psi(z)}{y_1(z)^2}\int_0^z|[y_1(\xi)y_2(\xi)]'|d\xi\nonumber\\
&=\frac{y_2(z)}{y_1(z)}\psi(z)+\frac{\psi(z)}{y_1(z)^2}(y_1(z)y_2(z)-1)\nonumber
\\ &\le 2\psi(z)\frac{y_2(z)}{y_1(z)}=2\psi(z)\beta(z),\quad
z\in\chi(x).\label{4.20}
\end{align}
Since for $z\in\chi(x)$ the function $\psi(x)$ satisfies the inequalities (see
\eqref{1.14})
\begin{align}
\psi(z)&=\sup_{t\in[0,z]}\left|\int_0^t(q(x+\xi)-q(x-\xi))d\xi\right|\le
\sup_{t\in\chi(x)}\left|\int_0^t(q(x+\xi)-q(x-\xi))d\xi\right|\nonumber\\
&\le \frac{F(x)}{\sqrt{k (x)}d(x)},\label{4.21}
\end{align}
by \eqref{4.21} we can continue estimate \eqref{4.20}:
\begin{align}
&|\beta'(z)| \le \frac{2F(x)}{\sqrt{k (x)}d(x)}\beta(z),\quad
z\in\chi(x), \quad x\in R\nonumber\\
&\Rightarrow\ -\frac{2F(x)}{\sqrt{k (x)}d(x)}\le \frac{\beta'(z)}{\beta(z)}\le
\frac{2F(x)}{\sqrt{k (x)}d(x)},\quad z\in\chi(x),\quad x\in R.\label{4.22}
\end{align}
Since $\beta(0)= 1,$ from \eqref{4.22} we get
\begin{equation}\label{4.23}
\exp(-2F(x))\le\beta(z)\le\exp(2F(x)),\quad z\in\chi(x),\quad x\in R.
\end{equation}

Let us check that $F(x)\to0$ as $|x|\to\infty.$ According to \eqref{1.14} and
\eqref{1.15}, we have \begin{align} &F(x) =\sqrt{k(x)} d(x)\sup_{z\in\chi(x)}
\left|\int_0^z(q(x+t)-q(x-t)dt\right|\nonumber\\
 &\le \frac{1}{\sqrt{k (x)}}\left[k (x)d(x)
\sup_{z\in[0,k (x)d(x)]}\left|\int_0^z(q(x+t)-q
(x-t))dt\right|\right]\le\frac{c_2}{\sqrt{k (x)}},\quad x\in R\label{4.24}
\end{align}
\begin{equation}\Rightarrow\ F(x)\to0\qquad\text{as}\quad |x|\to\infty.\label{4.25}
\end{equation}

Let now $|x|\ge t_0.$ Then from the assumption of the lemma and \eqref{4.24},
it follows that
\begin{equation}\label{4.26}
\alpha(x)\doe 2F(x)\le\frac{2c_2}{\sqrt{k (x)}}\le\frac{1}{4},\quad |x|\ge t_0.
\end{equation}

{}From \eqref{4.26}, we get
\begin{align}
&e^{\alpha x} =1+\alpha(x)+\sum_{n=2}^\infty\frac{(\alpha(x))^n}{n!}\le
1+\alpha(x)+\frac{\alpha(x)^2}{2}\sum_{k
=0}^\infty\left(\frac{\alpha(x)}{2}\right)^
k \nonumber\\
&=1+\alpha(x)+\frac{\alpha(x)^2}{2-\alpha(x)}\le
1+\alpha(x)+\frac{4}{7}(\alpha(x))^2\le 1+\alpha(x)+\frac{\alpha(x)}{7}\le
1+1.2\alpha(x).\label{4.27}
\end{align}
The lemma follows from \eqref{4.26}, \eqref{4.27} and \eqref{4.23}.
\end{proof}

\begin{lem}\label{lem4.7} Suppose $q(x)\in H,$ let $y(t)$ be the solution of
problem \eqref{4.3}--\eqref{4.4}, and let $t_0$ be the number from
\lemref{lem4.6}. Then for $|x|\ge t_0$, the following relations hold (see
\eqref{2.14}):
\begin{equation}\label{4.28}
\int_x^{\omega^+(x)}\frac{dt}{y(t)^2}=(1+\tau(x))\int_{\omega^-(x)}^x\frac{dt}{y(t)^2},\quad
|\tau(x)|\le 3.6F(x)\le\frac{3.6c_2}{\sqrt{k (x)}}.
\end{equation}
\end{lem}

\begin{proof}
By the definition of $t_0,$ for $|x|\ge t_0$ and $z\in\chi(x)$, we have the
following estimate for $|\gamma(z)|$ (see \eqref{4.18}):
\begin{equation}\label{4.29}
|\gamma(z)|\le 1.2F(x)\le\frac{1.2c_2}{\sqrt{k (x)}}\le
1.2\cdot\frac{1}{8}=0.15.
\end{equation}

{}From \eqref{4.18} and \eqref{4.29}, we get
\begin{align}
\int_x^{\omega^+(x)}&\frac{dt}{y(t)^2} =\int_0^{\sqrt{k
(x)}dx}\frac{dz}{y(x+z)^2}
=\int_0^{\sqrt{k (x)}d(x)}\frac{dz}{(1+\gamma(z))^2y(x-z)^2}\nonumber\\
&=\int_0^{\sqrt{ k (x)}d(x)}\frac{dz}{y(x-z)^2}-\int_0^{\sqrt{k
(x)}d(x)}\gamma(z)\frac
{2+\gamma(z)}{(1+\gamma(z))^2}\ \frac{dz}{y(x-z)^2}\nonumber\\
 &=\left[1-
 \int_0^{\sqrt{k (x)}d(x)}\gamma(z)\frac{2+\gamma(z)}{(1+\gamma(z))^2}\
\frac{dz}{y(x-z)^2}\left(\int_0^{\sqrt{k (x)}d(x)}\frac{dz}{y(x-z)^2}\right)
^{-1} \right]\nonumber\\
&\quad\cdot\int_0^{\sqrt{k (x)}d(x)}\frac{dz}{y(x-z)^2}\nonumber\\
&\doe(1+\tau(x))\int_0^{\sqrt{k (x)}d(x)}\frac{dt}{y(x-z)^2}=(1+\tau(x))\int_{
\omega^-(x)}^x\frac{dt}{y(t)^2},\quad |x|\ge t_0.\label{4.30}
\end{align}

It remains to prove the estimate $|\tau(x)|$ from \eqref{4.28}. We use
relations \eqref{4.29} and \eqref{4.18}:
$$|\tau(x)|\le\max_{z\in\chi(x)}|\gamma(x)|\frac{|2+\gamma(z)|}{(1+\gamma(z))^2}\le
1.2F(x)\frac{2.15}{0.85^2}\le 3.6F(x)\le\frac{3.6c_2}{\sqrt{k (x)}}.$$
\end{proof}

\begin{lem}\label{lem4.8} Suppose $q(x)\in H$, and let $t_0$ be the number from
\lemref{lem4.6}. Then for $|x|\ge t_0$, the following inequalities hold:
\begin{equation}\label{4.31}
|\rho'(x)|\le 3.6\left[F(x)+\sqrt{c_3}\exp\left(-\frac{\sqrt{k
(x)}}{3c_3}\right)\right]\le 11\frac{c_2+c_3\sqrt{c_3}}{k (x)},\quad |x|\ge
t_0. \end{equation}
\end{lem}

\begin{proof}
Below  when estimating $\varkappa(x)$, we use \eqref{4.11}, \eqref{4.6} and
\eqref{4.28}:
\begin{align}
\varkappa(x)&=\frac{v'(x)}{v(x)}\ \frac{u(x)}{|u'(x)|}=\frac{\tilde
u(x)}{\tilde
v(x)}(1+\nu(x))=\int_x^{\omega^+(x)}\frac{dt}{y(t)^2}\left(\int_{\omega^-(x)}^x
\frac{dt}{y(t)^2}\right)^{-1}(1+\nu(x))\nonumber\\
&=(1+\tau(x))(1+\nu(x))=1+\tau(x)+\nu(x)+\tau(x)\nu(x)\doe1+\mu(x),\quad |x|\ge
t_0.\label{4.32}
\end{align}
{}From \eqref{4.28} and \eqref{4.26}, we get
\begin{equation}\label{4.33}
|\tau(x)|\le 3.6\frac{c_2}{\sqrt{k(x)}}\le 0.45,\qquad |x|\ge t_0.
\end{equation}
Inequality \eqref{4.33}, together with \eqref{4.32}, \eqref{4.28} and
\eqref{4.11}, lead to the estimates:
\begin{align}
|\mu(x)|&\le |\tau(x)|+|\nu(x)|+|\tau(x)|\ |\nu(x)|\le |\tau(x)|+1.45|\nu(x)|
 \nonumber\\ &\le
3.6F(x)+3.6\sqrt{c_3}\exp\left(-\frac{\sqrt{k(x)}}{3c_3}\right)  \le
4F(x)+\frac{11c_3^{3/2}}{\sqrt{k(x)}}\nonumber\\ &\le
11\frac{c_2+c_3\sqrt{c_3}}{\sqrt{k(x)}},\quad |x|\ge t_0.\label{4.34}
\end{align}
Inequalities \eqref{4.31} follow from \eqref{4.34} and \eqref{4.1}.
\end{proof}

\begin{proof}[Proof of Theorems \ref{thm1.4} and \ref{thm1.5}]
Suppose $q(x)\in H.$ Then \eqref{4.31} implies \eqref{1.16} for $|x|\gg1.$ In
addition, \eqref{4.31} and \eqref{1.12} lead to \eqref{3.1}. Hence by
\lemref{lem3.1} we get \eqref{3.2}. We estimate $|\epsilon(x)|$ in \eqref{3.3}
using the estimate for $h(x).$ Below for $|x|\gg1$ we use \eqref{4.31}
\eqref{4.31}, \eqref{1.13} and \eqref{4.24}:
\begin{align}
h(x)&=\sup_{t\in\Delta(x)}|\rho'(t)|\le
3.6\sup_{t\in\Delta(x)}F(x)+3.6\sqrt{c_3}\sup_{t\in\Delta(x)}\left(-\frac{\sqrt{k(t)}}
{3c_3}\right)\nonumber\\
&\le
3.6\sup_{t\in\Delta(x)}F(t)+3.6\sqrt{c_3}\exp\left(-(3c_3c_1^{1/2})^{-1}\sqrt{k(x)}\right)
\nonumber\\
&\le
c\left\{\sup_{t\in\Delta(x)}F(t)+\exp\left(-c^{-1}\sqrt{k(x)}\right)\right\}=c\beta(x).
\label{4.35}
\end{align}

{}From \eqref{4.35} and \eqref{3.3}, we get \eqref{1.19}.  {}From \eqref{4.24}
and \eqref{1.13}, we obtain
\begin{align*}
\beta(x)&=\sup_{t\in\Delta(x)}F(t)+\exp\left(-c^{-1}\sqrt{k(x)}\right)\le\sup_{t\in\Delta(x)}
\frac{c_2}{\sqrt{k(t)}}+\frac{c}{\sqrt{k(x)}}\\
&\le
\frac{c_2\sqrt{c_1}}{\sqrt{k(x)}}+\frac{c}{\sqrt{k(x)}}\le\frac{c}{\sqrt{k(x)}}
\quad\Rightarrow\quad\eqref{1.20}.
\end{align*}
Thus \thmref{thm1.4} is proved.

To prove \thmref{thm1.5}, we set $|x|\ge s_1>s_0$ (see \eqref{1.26} --
\eqref{1.27}). Then $s_0>t_0$ because of \eqref{1.26}, where $t_0$ is the
number from \lemref{lem4.6}. Hence $|\rho'(x)|\le 10^{-3}$ according to
\eqref{4.31} and \eqref{1.26}. Formula \eqref{3.2} is proved similarly to
\thmref{thm1.4}. Since \eqref{4.31} coincides with \eqref{1.27}, it remains to
estimate $|\varepsilon(x)|$ using \eqref{3.6}. We use one of the inequalities
\eqref{4.35} and obtain:
$$|\varepsilon(x)|\le 18h(x)\le
65\left(\sup_{t\in\Delta(x)}F(t)+\sqrt{c_3}\exp\left(-\left(3c_3\sqrt{c_1}\right)^{-1}\sqrt{
k(x)}\right)\right)=\eta_2(x).$$ \thmref{thm1.5} is proved.
\end{proof}

\section{Comparison of two asymptotic formulas}

In this section we prove \thmref{thm1.7}.

\begin{proof}[Proof of \thmref{thm1.7}]
Suppose $q(x)\in H$ and $q(x)\ge 1$ for $x\in R.$ Consider \eqref{1.17} and
\eqref{1.19}.  In these relations $\alpha(x)$ and $\beta(x)$ are positive,
continuous for $x\in R$ functions, $\beta(x)\le\alpha(x)$ for $x\in R$ and
$\alpha(x)\to0,$\ $\beta(x)\to0$ as $|x|\to\infty.$ Therefore the theorem will
be proved if (see   \eqref{1.31})
\begin{equation}\label{5.1}
L=\sup_{q(\cdot)\in
H}\varlimsup_{|x|\to\infty}\frac{\alpha(x)}{\beta(x)}=\infty.
\end{equation}
Note that to prove \eqref{5.1}, it suffices to give an example of the function
$q(\cdot)$ which, on the one hand, satisfies the above-mentioned assumptions
and, on the other hand, for the functions $\alpha(x)$ and $\beta(x)$
constructed by $q(\cdot)$ the following equality holds:
\begin{equation}\label{5.2} \tilde
L=\varlimsup_{|x|\to\infty}\frac{\alpha(x)}{\beta(x)}=\infty.
\end{equation}
(Indeed, if \eqref{5.2} holds, then $\infty=\tilde L\le L\le \infty\
\Rightarrow\ \tilde L=L=\infty.)$

Let us construct such a function. Denote
\begin{equation}\label{5.3}
\sigma_n=\left[\left.\sigma_n^{(-)},\sigma_n^{(+)}\right)\right.=[n^2,(n+1)^2),\quad
q_n=\left(1+\frac{1}{n}\right)^n,\quad n=1,2,\dots
\end{equation}
Suppose $q(-x)=q(x)$ for $x\ge0$ and
\begin{equation}\label{5.4}
q(x)=\begin{cases} q_n\quad &\text{if}\quad x\in\sigma_n,\quad n=1,2,\dots\\
2\quad &\text{if}\quad x\in[0,1).
\end{cases}
\end{equation}

Clearly, we have $1\le q(x)\in L_1^{\loc}(R),$\ $x\in R.$ We show that $q(x)\in
H.$ We need to estimate the function $d(x).$ Since
\begin{equation}\label{5.5}
2\le q(x)\le3 \qquad\text{for}\quad x\in R,
\end{equation}
by \eqref{1.11} we have
$$2=d(x)\int_{x-d(x)}^{x+d(x)}q(t)dt\ge d(x)\int_{x-d(x)}^{x+d(x)}2dt=4d(x)^2\
\Rightarrow\ d(x)\le\frac{1}{\sqrt{2}},\ x\in R,$$
$$2=d(x)\int_{x-d(x)}^{x+d(x)}q(t)dt\le d(x)\int_{x-d(x)}^{x+d(x)}3dt=6d(x)\
\Rightarrow\ d(x)\le\frac{1}{\sqrt{3}},\ x\in R.$$ Hence
\begin{equation}\label{5.6}
\frac{1}{\sqrt{3}}\le d(x)\le\frac{1}{\sqrt{2}},\qquad x\in R.
\end{equation}

We introduce the function
\begin{equation}\label{5.7}
k(x)=\begin{cases} \sqrt{|x|},\quad &\text{if}\quad |x|\ge 4\\
 2,\quad &\text{if}\quad |x|\le 4\end{cases}
 \end{equation}
 Let us check that in case \eqref{5.7} all the assumptions of Definition
 \ref{defn1.2} are satisfied. {}From \eqref{5.6} and \eqref{5.7}, we get
 relations \eqref{1.12} and \eqref{1.13}. In particular, \eqref{1.12}
 immediately follows from \eqref{5.7}. We prove \eqref{1.13}. Let $x\ge9.$ Then
 \begin{gather*}
 x-\sqrt{x}=x\left(1-\frac{1}{\sqrt{x}}\right)\ge9\left(1-\frac{1}{3}\right)=6\ge
 4\\
 \Rightarrow\quad
 [x-k(x)d(x),x+k(x)d(x)]\subseteq\left[x-\sqrt{x},x+\sqrt{x}\right],\\
 \left[x-\sqrt{x},x+\sqrt{x}\right]\cap[-4,4]=\emptyset.
 \end{gather*}
 Thus for $x\ge9,$ inequalities \eqref{1.13} are true.

  The estimates proved for $|x|\ge9$ can be easily extended to the whole number
  axis using \lemref{lem2.6}. It remains to check \eqref{1.14}. Consider
  $\Phi(x)$ (see \eqref{1.14}) for $x\in\sigma_n,$\ $n\ge 2.$
  Clearly, if $x\in\sigma_n,$ then
$$\left.\begin{array}{ll}
  x+k(x)d(x)\le x+\sqrt x\le (n+1)^2+n+1<(n+2)^2\nonumber\\
  x-k(x)d(x)\ge x-\sqrt x\ge n^2-n>(n-1)^2;\nonumber\\
  \end{array}\right\}\ \Rightarrow$$
\begin{equation}\label{5.8}
  [x-k(x)d(x),x+k(x)d(x)]\subseteq\sigma_{n-1}\cup\sigma_n\cup\sigma_{n+1}\quad\text{for}
  \quad x\in\sigma_n,\quad n\ge 2.
  \end{equation}
  Now from the condition $x\in\sigma_n$ and \eqref{5.8}, \eqref{5.4},
  \eqref{5.7}, \eqref{5.6} and \eqref{1.14}, we get
  \begin{align}
  \Phi(x)&=k(x)d(x)\sup_{x\in[0,k(x)d(x)]}\left|\int_0^z(q(x+t)-q(x-t)dt\right|
\nonumber\\
 &\le\frac{n+1}{\sqrt{2}}\max\left\{|q_{n+1}-q_n|,|q_n-q_{n-1}|\right\}\cdot\frac{n+1}{\sqrt{2}}
\nonumber\\
&=\frac{(n+1)^2}{2}\max\left\{|q_{n+1}-q_n|,|q_n-q_{n-1}|\right\}.\label{5.9}
\end{align}

Since the following inequalities hold (see \cite[Section I, problem 170]{11}):
\begin{equation}\label{5.10}
\left.\begin{array}{ll}
\qquad \frac{e}{2n+2}<e-q_n<\frac{e}{2n+1} \\
-\frac{e}{2n+3}<q_{n+1}-e<-\frac{e}{2n+4}
\end{array}\right\}\
\Rightarrow \frac{e}{(2n+2)(2n+3)}<q_{n+1}-q_n<\frac{3e}{(2n+1)(2n+4)}
\end{equation}
by \eqref{5.9} and \eqref{5.10}, we obtain
\begin{equation}\label{5.11}
\Phi(x)\le\frac{(n+1)^2}{2}\ \frac{3e}{(2n-1)(2n+2)}\le
C<\infty\qquad\text{for}\quad x\in\sigma_n,\quad n\ge 2.
\end{equation}

We omit the obvious proof of \eqref{1.14} using \eqref{5.11}. Thus $q(x)\in H,$
and it remains to prove \eqref{5.2}. Let (see \eqref{5.3})
\begin{equation}\label{5.12}
x_n=\frac{\sigma_n^{(-)}+\sigma_n^{(+)}}{2}=n^2+n+\frac{1}{2},\qquad n\ge 1.
\end{equation}
Let us compute $\sup\limits_{t\in\Delta(x_n)}F(t).$ Note that if
$t\in\Delta(x_n),$ then \eqref{5.6} implies elementary inequalities
$$\left.\begin{array}{ll}
  t+\sqrt{k(x)}d(t)\le x_n+\frac{1}{\sqrt2}+\sqrt{k\left(x_n+\frac{1}{\sqrt2}\right)}
  \frac{1}{\sqrt2}<(n+1)^2\quad\text{for}\quad n\gg1,\nonumber\\
 t-\sqrt{k(x)}d(t)\ge x_n-\frac{1}{\sqrt2}-\sqrt{k\left(x_n-\frac{1}{\sqrt2}\right)}
  \frac{1}{\sqrt2}>n^2\quad\text{for}\quad n\gg1,\nonumber\\
  \end{array}\right\}\ \Rightarrow$$
\begin{equation}\label{5.13}
 \left[t-\sqrt{k(x)}d(t),t+\sqrt{k(x)}d(t)\right]\subset\sigma_n\quad\text{for}\quad t\in\Delta_n
 \quad\text{and}\quad n\gg 1.
  \end{equation}
  Furthermore, according to \eqref{5.13} and \eqref{5.4}, we have
  $$q(t+\xi)=q(t-\xi)=q_n\qquad\text{for}\qquad |\xi|\le\sqrt{k(t)}d(t),\quad
  t\in\Delta(x_n),\quad n\gg1,$$
  and therefore $\sup\limits_{t\in\Delta(x_n)}F(t)=0$ (see \eqref{1.15}). By
  \eqref{1.20}, this implies
  \begin{equation}\label{5.14}
  \beta(x_n)=\exp\left(-c^{-1}\sqrt{k(x_n)}\right),\qquad n\gg1.
  \end{equation}

  Now we consider the value
  $$\tilde
  F(s)\bigm|_{s=\sigma_n^{(+)}}=\sqrt{k(x)}d(s) \left|\left.\int_0^{\sqrt{k(s)}ds}(q(s+t)-q(s-t))
  dt\right|\
  \right|_{s=\sigma_n^{(+)}}.
 $$
 Since for
 $t\in\left(\left.0,\sqrt{k\left(\sigma_n^{(+)}\right)}d\left(\sigma_n^{(+)}\right)\right]
 \right. $, we have
 $$q(\sigma_n^{(+)}+t)=q_{n+1},\qquad q(\sigma_n^{(+)}-t))=q_n,$$
 using \eqref{5.10}, \eqref{5.6} and \eqref{5.7}, we get
 $$\tilde
 F(\sigma_n^{(+)})=k(\sigma_n^{(+)})d(\sigma_n^{(+)})^2(q_{n+1}-q_n)\ge\frac{n+1}{3}\
 \frac{e}{(2n+2)(2n+3)}=\frac{e}{6(2n+3)}.$$
 The last inequality yields the estimates
 \begin{equation}\label{5.15}
 \sup_{t\ge x_n-d(x_n)}F(t)\ge F(\sigma_n^{(+)})\ge\tilde F(\sigma_n^{(+)})\ge
 \frac{e}{6(2n+3)},\qquad n\gg1.
 \end{equation}
 {}From \eqref{5.15} and \eqref{1.18}, we finally get
 \begin{equation}\label{5.16}
 \alpha(x_n)\ge\exp\left(-c^{-1}\sqrt{k(x_n)}\right)+\frac{e}{6(2n+3)},\qquad n\gg1.
 \end{equation}
 Relations \eqref{5.14} and \eqref{5.15} imply \eqref{5.1}.
 Indeed,
 $$\frac{\alpha(x_n)}{\beta(x_n)}\ge
1+\frac{e}{6(2n+3)}\exp\left(c^{-1}\sqrt{k(x_n)}\right)\qquad\text{for}\quad
 n\gg1$$
 $$\Rightarrow\quad\lim_{n\to\infty}\frac{\alpha(x_n)}{\beta(x_n)}=\infty\quad
 \Rightarrow\quad \tilde L=\infty\quad\Rightarrow\quad L=\infty.$$
\end{proof}

\section{Properties of solutions of the Riccati equation}

In this section, we prove \thmref{thm1.8}. Below we use the following
assertion.

\begin{thm}\label{thm6.1} \cite[\S402]{9}
The general solution of equation \eqref{1.32} is of the form
\begin{equation}\label{6.1}
y(x)=\frac{c_1v'(x)+c_2u'(x)}{c_1v(x)+c_2u(x)}.
\end{equation}\end{thm}
Here $\{u(x),v(x)\}$ is a PFSS of equation \eqref{1.1}, $c_1,c_2$ are arbitrary
constant, $|c_1|+|c_2|\ne0$.

 \begin{proof}[Proof of \thmref{thm1.8}]
 Let $q(x)\in H.$ Set
 \begin{equation}\label{6.2}
 y_2(x)=\frac{v'(x)}{v(x)},\qquad y_1(x)=\frac{u'(x)}{u(x)},\qquad x\in R.
 \end{equation}
 Then by \eqref{3.8} and \thmref{thm1.4}, we get \eqref{1.33}:
 $$\lim_{|x|\to\iy}y_2(x)d(x)=\lim_{|x|\to\iy}\frac{v'(x)}{v(x)}d(x)=\lim_{|x|\to\iy}(1+
 \rho'(x))\frac{d(x)}{2\rho(x)}=1,$$
 $$\lim_{|x|\to\iy}y_1(x)d(x)=\lim_{|x|\to\iy}\frac{u'(x)}{u(x)}d(x)=\lim_{|x|\to\iy}(
 \rho'(x)-1)\frac{d(x)}{2\rho(x)}=-1.$$

Consider the second part of assertion A). Suppose that there exists a solution
$y(x)$ of equation \eqref{1.33} which satisfies the following properties of the
solution $y_2(x):$
\begin{enumerate}\item[1)] the solution $y(x)$ is defined for all $x\in R;$
\item[2)] the following equalities hold:
\begin{equation}\label{6.3}
\lim_{x\to-\iy}y(x)d(x)=\lim\limits_{x\to\iy}y(x)d(x)=1.
\end{equation}
\end{enumerate}
Let us show that 1) and 2) imply $y(x)\equiv y_2(x)$ for $x\in R.$

We need the following assertion.
\begin{lem}\label{lem6.2}
Suppose that conditions \eqref{1.2} hold, and let $y(x)$ be a solution of
equation \eqref{1.32} such that $y(x)\ne y_1(x),$\ $y(x)\ne y_2(x).$ Then if
the solution $y(x)$ is defined for all $x\in R,$ then
\begin{equation}\label{6.4}
y_1(x)<y(x)<y_2(x)\qquad\text{for}\qquad x\in R.\end{equation}
\end{lem}

\begin{proof}
Suppose that $y(x_0)>y_2(x_0)$ for some $x_0\in R.$ By the hypothesis of the
lemma, in representation \eqref{6.1} we have $c_1\ne0,$\ $c_2\ne0,$ and
therefore
\begin{equation}\label{6.5}
y(x)=\frac{v'(x)+\theta u'(x)}{v(x)+\theta u(x)},\qquad \theta\ne0,\quad
\theta=\frac{c_2}{c_1},\quad x\in R.
\end{equation}
Since $y(x_0)>y_2(x_0),$ \eqref{1.3} and \eqref{1.4} imply
\begin{align*}
0&<y(x_0)-y_2(x_0)=\frac{v'(x_0)+\theta u'(x_0)}{v(x_0)+\theta
(x_0)}-\frac{v'(x_0)}{v(x_0)}=-\frac{\theta}{(v(x_0+\theta u(x_0))v(x_0)}\\
&=-\frac{1}{\theta^{-1}+u(x_0)/v(x_0)}\
\frac{1}{v(x_0)^{2}}\quad\Rightarrow\quad
\theta^{-1}+\frac{u(x_0)}{v(x_0)}<0\quad\Rightarrow\quad \theta<0.
\end{align*}

Let $\varphi(x)=u(x)/v(x)$,\ $x\in R.$ According to \eqref{1.3}, \eqref{1.4}
and \eqref{1.5}, this function satisfies the properties $\varphi(x)\to\iy$ as
$x\to-\iy,$\ $\varphi(x)\to0$ as $x\to\iy$
$$\varphi'(x)=\frac{u'(x)v(x)-v'(x)u(x)}{v(x)^2}=-\frac{1}{v(x)^2}<0\qquad\text{for}\quad
x\in R.$$ Hence there exists $x_1$ such that $\varphi(x_1)=-\theta^{-1},$ or,
equivalently, \begin{equation}\label{6.6} v(x_1)+\theta(x_1)=0.
\end{equation}
Together with equality \eqref{6.6}, the following inequality holds:
\begin{equation}\label{6.7}
v'(x_1)+\theta u'(x_1)=v'(x_1)+| \theta u'(x_1)|>0
\end{equation}
(see \eqref{3.9}). {}From \eqref{6.6}, \eqref{6.7} and \eqref{6.5}, it follows
that the solution $y(x)$ is not defined for $x=x_1; $ contradiction. Hence
$y(x)\le y_2(x)$ for all $x\in R.$ But $y(x)\ne y_2(x)$ by hypothesis which
leads to the upper estimate in \eqref{6.4}. The second inequality of
\eqref{6.4} can be checked in a similar way.
 \end{proof}

 \begin{cor}\label{cor6.3} Assuming the hypothesis of \lemref{lem6.2}, the solution
  of equation \eqref{1.32} is of the form \eqref{6.5} with
 $\theta>0.$
 \end{cor}

\begin{proof} Taking into account all that was mentioned above, it only remains
to check that $\theta>0.$ {}From \eqref{6.4}, \eqref{1.3} and \eqref{1.4}, it
follows that
$$
0 <y(x)-y_1(x)=\frac{v'(x)+\theta u'(x)}{v(x)+\theta
u(x)}-\frac{u'(x)}{u(x)}=\frac{1}{u(x)(v(x)+\theta u(x))},\qquad x\in R;$$
$$
0 <y_2(x)-y(x)=\frac{v'(x)}{v(x)}-\frac{v'(x)+\theta u'(x)}{v(x)+\theta
u(x)}=\frac{\theta}{(v(x)(v(x)+\theta u(x))}.$$ The first inequality implies
$v(x)+\theta u(x)>0,$\ $ x\in R.$ Then $\theta>0$ in view of the second
inequality.
\end{proof}

We can now finish the proof of assertion A). First note that if $q(x)\in H,$
then in addition to \eqref{1.5} we have the relations
\begin{equation}\label{6.8}
\lim_{x\to-\infty}\frac{v'(x)}{u'(x)}=\lim_{x\to\iy}\frac{u'(x)}{v'(x)}=0.
\end{equation}
Indeed, from \eqref{3.9}, \eqref{3.8}, \eqref{1.16}, \eqref{1.12} and
\eqref{1.5}, it follows that
\begin{align*}
\lim_{x\to-\iy}\frac{v'(x)}{u'(x)}&=\lim_{x\to-\iy}\frac{v'(x)}{v(x)}\
\frac{v(x)}{u(x)}\
\frac{u(x)}{u'(x)}=\lim_{x\to-\iy}\frac{1+\rho'(x)}{2\rho(x)}\
\frac{v(x)}{u(x)}\ \frac{2\rho(x)}{\rho'(x)-1}\\
&=\lim_{x\to-\iy}\frac{1+\rho'(x)}{\rho'(x)-1}\lim_{x\to-\iy}\frac{v(x)}{u(x)}=0.
\end{align*}

The second equality of \eqref{6.8} can be proved in a similar way. Let $y(x)$
be a solution of \eqref{1.32} which does not coincide with $y_2(x)$ for $x\in
R$ and satisfies properties 1)--2) (see above). Then by \corref{cor6.3} the
solution $y(x)$ is of the form \eqref{6.5} with $\theta>0.$ In the following
relations, we use \eqref{6.3}, \eqref{6.5}, \eqref{3.9}, \eqref{1.4},
\eqref{6.8}, \eqref{1.5} and
\begin{align*}
1&=\lim_{x\to-\iy}y(x)d(x)=\lim_{x\to-\iy}\frac{v'(x)+\theta u'(x)}{v(x)+\theta
u(x)}
d(x)=\lim_{x\to-\iy}\frac{\theta^{-1}\frac{v'(x)}{u'(x)}+1}{\theta^{-1}\frac{v(x)}{u(x)}+1}
\ \frac{u'(x)}{u(x)}d(x)\\
&=\lim_{x\to-\iy}y_1(x)d(x)=-1.
\end{align*}
Contradiction.  Hence $y(x)=y_2(x),$\ $x\in R.$ The part of assertion A)
related to $y_1(x)$ can be proved similarly.

Let us prove B). Note that $y_+(x)=y_1(x)$ if and only if $c_1=0$ in
\eqref{6.1}. In fact, for $x\in [c,\iy)$ we have
\begin{align*}
&y_+(x)
=y_1(x)\quad\Leftrightarrow\quad\frac{c_1v'(x)+c_2u'(x)}{c_1v(x)+c_2u(x)}=\frac{u'(x)}
{u(x)}\\
&\Leftrightarrow\quad c_1(v'(x)u(x)-u'(x)v(x))=0\quad\Leftrightarrow\quad
c_1=0.
\end{align*}
Thus the condition $y_+(x)\ne y_1(x)$ implies that in this case we have
$c_1\ne0$ in \eqref{6.1}. Hence, as in the proof of A) given above, we obtain
\begin{align*}
\lim_{x\to\iy}y_+(x)d(x)&=\lim_{x\to\iy}\frac{c_1v'(x)+c_2u'(x)}{c_1v(x)+c_2u(x)}d(x)\\
&=\lim_{x\to\iy}\frac{1+\frac{c_2}{c_1}\
\frac{u'(x)}{v'(x)}}{1+\frac{c_2}{c_1}\ \frac{u(x)}{v(x)}}\
\frac{v'(x)}{v(x)}dx=\lim_{x\to\iy}y_2(x)d(x)=1.
\end{align*}

The converse statement is an obvious consequence of \eqref{1.34}. Assertion C)
can be proved in the same way as assertion B).
 \end{proof}

 \section{Asymptotics of the Otelbaev function at infinity}

 The problem that is considered in this section arises as a result of attempts
 to use \linebreak Theorems~\ref{thm1.4} and \ref{thm1.5} in order to study concrete
 equations \eqref{1.1} and \eqref{1.32}. It is easily seen that to study
 theoretical problems related to asymptotic behaviour of the function $\rho(x)$
 at infinity, one can use formula \eqref{1.19} without additional restrictions
 to $q(x)$. (See, for example, \thmref{thm1.8}. Another such example was given
 in \cite{6} where \eqref{1.17} helped to find asymptotics at infinity of the
 distribution function of the spectrum of the Sturm-Liouville operator.)
 However, to apply Theorems \ref{thm1.4} and \ref{thm1.5} to concrete
 equations, one has to know the asymptotic estimates of $d(x)$ for
 $|x|\to\infty.$
 The proof of such estimates is a separate technical problem which is not at
 all related to the initial question on the properties of $\rho(x)$ for
 $|x|\to\infty.$ To solve this problem, additional requirements different from
 the conditions of \thmref{thm1.4} are imposed on the function $q(x).$ In
 \cite{6}, such a requirement is condition \eqref{2.1}.

 In the following theorem, we find an asymptotics of $d(x)$ at infinity under
 condition \eqref{1.2} and some additional requirements which are more
 convenient for   practical checking than the corresponding conditions from
 \cite{6}.

 \begin{thm}\label{thm7.1} Suppose that $0\le q(x)\in L_1^{\loc}(R),$\ $x\in R$
 and one can represent the function $q(x)$ in the form
 \begin{equation}\label{7.1}
 q(x)=q_1(x)+q_2(x),\qquad x\in R,
 \end{equation}
 where $q_1(x)$ is positive for $x\in R$ and twice differentiable for $|x|\gg1$
 and $q_2(x)\in L_1^{\loc}(R).$
 Denote
 \begin{equation}\label{7.2}
 A(x)=[0,2q_1(x)^{-1/2}],\quad x\in R,
 \end{equation}
 \begin{equation}\label{7.3}
 \varkappa_1(x)=\frac{1}{q_1(x)^{3/2}}\sup_{t\in
 A(x)}\left|\int_{x-t}^{x+t}q_1''(\xi)d\xi\right|,\quad |x|\gg1,
 \end{equation}
 \begin{equation}\label{7.4}
 \varkappa_2(x)=\frac{1}{\sqrt{q_1(x)}}\sup_{t\in
 A(x)}\left|\int_{x-t}^{x+t}q_2(\xi)d\xi\right|,\quad x\in R.
 \end{equation}
Then if the following condition holds:
\begin{equation}\label{7.5}
\varkappa_1(x)\to0,\qquad \varkappa_2(x)\to0\qquad\text{as}\qquad |x|\to\iy,
\end{equation}
condition \eqref{1.2} also holds and for every $x\in R$, the equation
\eqref{1.11} has a unique positive solution $d(x).$ Moreover,
\begin{equation}\label{7.6}
d(x)=\frac{1+\delta(x)}{\sqrt{q_1(x)}},\quad |\delta(x)|\le
2(\varkappa_1(x)+\varkappa_2(x))\quad\text{for}\quad |x|\gg 1,
\end{equation}
\begin{equation}\label{7.7}
c^{-1}\le d(x)\sqrt{q_1(x)}\le c\qquad\text{for}\qquad x\in R.
\end{equation}
 \end{thm}

 \begin{proof}[Proof of \thmref{thm7.1}]
\quad{}

 We need the following lemma.

 \begin{lem}\label{lem7.2}
 Suppose that $q_1$ satisfies the hypothesis of \thmref{thm7.1}. For a given
 $x\in R,$ consider the following equation in $\hat d\ge0:$
 \begin{equation}\label{7.8}
 S(\hat d)=2,\qquad S(\hat d)=\hat d\int_{x-\hat d}^{x+\hat d} q_1(t)dt.
 \end{equation}
Equation \eqref{7.8} has a unique positive solution $\hat d(x);$ moreover (see
\eqref{7.3})
\begin{equation}\label{7.9}
\hat d(x)=\frac{1+\delta_1(x)}{\sqrt{q_1(x)}},\qquad
|\delta_1(x)|\le\varkappa_1(x)\qquad \text{for}\qquad |x|\gg 1,
\end{equation}
\begin{equation}\label{7.10}
c^{-1}\le \hat d(x)\sqrt{q_1(x)}\le c\qquad\text{for}\qquad x\in R.
\end{equation}
\end{lem}

\begin{proof} Clearly, $S(0)=0,$\ $S(\hat d)\to\iy$ as $\hat d\to\iy$ and $S(\hat
d)$ is monotone increasing in $\hat d\ge0.$ This implies that for every $x\in
R$ equation \eqref{7.8} has a unique positive solution $\hat d(x).$ To estimate
$\hat d(x),$ let us write the function $S(\hat d)$ in the form \eqref{7.11}:
\begin{align}
S(\hat d)&=\hat d\int_{x-\hat d}^{x+\hat d} q_1(t)dt=\hat d\int_0^{\hat
d}[q_1(x+t)+q_1(x-t)]dt \nonumber\\ &=2q_1(x)\hat d^2+\hat d\int_0^{\hat
d}[q_1(x+t)-2q_1(x)+q_1(x-t)]dt  \nonumber\\ &=2q_1(x)\hat d^2+\hat
d\int_0^{\hat d}\int_0^t\int_{x-\xi}^{x+\xi}q_1''(s)dsd\xi dt.\label{7.11}
\end{align}

Set (see \eqref{7.3})
$$\eta(x)=(1+\varkappa_1(x))q_1(x)^{-1/2},\qquad |x|\gg1.$$
By \eqref{7.5}, $\varkappa_1(x)\le 1$ for all $|x|\gg1,$ and therefore $\eta
(x)\in A(x)$ for all $|x|\gg1$ (see \eqref{7.2}). Then from \eqref{7.11} and
\eqref{7.3}, it follows that
\begin{align}
S(\eta(x))&=2(1+\varkappa_1(x))^2+\frac{1+\varkappa_1(x)}{\sqrt{q_1(x)}}\int_0^{\eta(x)}\int
_0^t\int_{x-\xi}^{x+\xi}q_1''(s)dsd\xi dt\nonumber\\
&\ge 2(1+\varkappa_1(x))^2-\frac{1}{2}\
\frac{(1+\varkappa_1(x))^2}{q_1(x)^{3/2}}\sup_{\xi\in
A(x)}\left|\int_{x-\xi}^{x+\xi}q_1''(s)ds\right|\nonumber\\
&=2(1+\varkappa_1(x))^2-\frac{\varkappa_1(x)(1+\varkappa_1(x))^3}{2}\ge
2+2\varkappa_1(x)\ge 2.\label{7.12}
\end{align}
 By \lemref{lem2.2}, \eqref{7.12} implies the inequality
\begin{equation}\label{7.13}
\hat d(x)\le\eta(x)=(1+\varkappa_1(x))q_1(x)^{-1/2}\qquad\text{for all}\qquad
|x|\gg1.
\end{equation}

Let now
$$\eta(x)=(1+\varkappa_1(x))^{-1}q_1(x)^{-1/2},\qquad |x|\gg 1.$$
Then, as above, we have $\eta(x)\in A(x)$ and using \eqref{7.11} and
\eqref{7.3}, we obtain
\begin{align}
S(\eta(x))&=\frac{2}{(1+\varkappa_1(x)^2}+\frac{1}{1+\varkappa_1(x)}\
\frac{1}{\sqrt{q_1(x)}}\int_0^{\eta(x)}\int_0^t\int_{x-\xi}^{x+\xi}q_1''(s)dsdtd\xi\nonumber\\
&\le \frac{2}{1+\varkappa_1(x))^2}+\frac{1}{2}\
\frac{\eta(x)^2}{1+\varkappa_1(x)}\ \frac{1}{\sqrt{q_1(x)}}\sup_{\xi\in
A(x)}\left|\int_{x-\xi}^{x+\xi}q_1''(s)ds\right|\nonumber\\
&\le\frac{2}{(1+\varkappa_1(x))^2}+\frac{\varkappa_1(x)}{2(1+\varkappa_1(x))^3}\le
2.\label{7.14}
\end{align}
By \lemref{lem2.2}, \eqref{7.14} implies the inequality
\begin{equation}\label{7.15}
\hat d(x)\ge \eta(x)=(1+\varkappa_1(x)^{-1}q_1(x)^{-1/2}\qquad\text{for
all}\qquad |x|\gg1.
\end{equation}
Estimates \eqref{7.13} and \eqref{7.15} yield \eqref{7.9}. Inequalities
\eqref{7.10} follows from \eqref{7.9} and \lemref{lem2.6}.
\end{proof}

We now prove \thmref{thm7.1}.  Consider the following equation in $d\ge0:$
\begin{equation}\label{7.16}
S(d)=2,\qquad S(d)=d\int_{x-d}^{x+d}q(\xi)d\xi.
\end{equation}
Let $\eta(x)=(1+\varkappa_2(x)\hat d(x)$ and $|x|\gg1$ (see \eqref{7.4}).
{}From \eqref{7.5} and \eqref{7.9}, it follows that $\eta(x)\in A(x)$ for all
$|x|\gg1.$ For such an $x$, \lemref{lem7.2} and \lemref{7.4} imply
\begin{align}
S(\eta(x))&=\eta(x)\int_{x-\eta(x)}^{x+\eta(x)}q(t)dt=\eta(x)\int_{x-\eta(x)}^{x+\eta(x)}q_1(t)dt
+\eta(x)\int_{x-\eta(x)}^{x+\eta(x)}q_2(t)dt\nonumber\\
&\ge (1+\varkappa_2(x))\hat d(x)\int_{x-\hat d(x)}^{x+\hat d(x)}
q_1(t)dt+\eta(x)\int_{x-\eta(x)}^{x+\eta(x)}q_2(t)dt\nonumber\\
&\ge
2(1+\varkappa_2(x))-\frac{(1+\varkappa_1(x))(1+\varkappa_2(x))}{\sqrt{q_1(x)}}\left|\int_{x-
\eta(x)}^{x+\eta(x)}q_2(t)dt\right|\nonumber\\ &\ge
2(1+\varkappa_2(x))-\varkappa_2(x)(1+\varkappa_2(x))(1+\varkappa_1(x))\ge
2.\label{7.17}
\end{align}

{}From \eqref{7.17} and the definition \eqref{7.16} of the function $S(d)$, it
is not hard to conclude (see \S2) that for all $|x|\gg1$ equation \eqref{7.16}
has a unique positive root $d(x).$ This implies that this property of equation
\eqref{7.16} remains true for all $x\in R$ and therefore, in particular,
\eqref{1.2} holds. Furthermore, \eqref{7.17} and \lemref{lem2.2} lead to the
inequality
\begin{equation}\label{7.18}
d(x)\le\eta(x)=(1+\varkappa_2(x))\hat d(x)\qquad\text{for all}\qquad |x|\gg1.
\end{equation}
Set
$$\eta(x)=(1+\varkappa_2(x))^{-1}\hat d(x)\qquad\text{for}\qquad |x|\gg1.$$
{}From \eqref{7.5} and \eqref{7.9}, it follows that $\eta(x)\in A(x)$ for all
$|x|\gg1.$ Then according to \eqref{7.4}, \lemref{lem7.2} implies
\begin{align}
S(\eta(x))&=\eta(x)\int_{x-\eta(x)}^{x+\eta(x)}q(t)dt=\eta(x)\int_{x-\eta(x)}^{x+\eta(x)}q_1(t)
dt+\eta(x)\int_{x-\eta(x)}^{x+\eta(x)}q_2(t)dt\nonumber\\
&\le\frac{1}{1+\varkappa_2(x)}\hat d(x)\int_{x-\hat d(x)}^{x+\hat
d(x)}q_1(t)dt+\eta(x)\int_{x-\eta(x)}^{x+\eta(x)}q_2(t)dt\nonumber\\
&\le\frac{2}{1+\varkappa_2(x)}+\frac{1+\varkappa_1(x)}{1+\varkappa_2(x)}\
\frac{1}{\sqrt{q_1(x)}}\sup_{t\in
A(x)}\left|\int_{x-t}^{x+t}q_2(\xi)d\xi\right|\nonumber\\
&=\frac{2}{1+\varkappa_2(x)}+\varkappa_2(x)\frac{1+\varkappa_1(x)}{1+\varkappa_2(x)}\le
2.\label{7.19}
\end{align}
Hence by \lemref{lem2.2} and estimate \eqref{7.19}, we have
\begin{equation}\label{7.20}
d(x)\ge\eta(x)=(1+\varkappa_2(x))^{-1}\hat d(x)\qquad\text{for all}\qquad
|x|\gg1.
\end{equation}

Set
$$d(x) =(1+\alpha(x))\hat d(x),\qquad|x|\gg1.$$
Then using the facts proved above, we obtain $|\alpha(x)|\le\varkappa_2(x)$ for
all $|x|\gg1.$ Therefore, taking into account \eqref{7.9}, for all $|x|\gg1, $
we get
\begin{gather*}
 d(x) =(1+\alpha(x))\hat
d(x)=\frac{(1+\alpha(x))(1+\delta_1(x))}{\sqrt{q_1(x)}}:=\frac{1+\delta(x)}{\sqrt{q_1(x)}}\\
 \Rightarrow\ |\delta(x)|\le|\alpha(x)|+|\delta_1(x)|+|\alpha(x)\delta_1(x)|\le
2(|\alpha(x)|+|\delta_1(x)|)\le 2(\varkappa_1(x)+\varkappa_2(x))\ \Rightarrow\
\eqref{7.6}.
\end{gather*}
Inequalities \eqref{7.7} follows from \eqref{7.6} and \lemref{lem2.6}.
\end{proof}

\section{Example}

In this section, we consider equation \eqref{1.1} and \eqref{1.33} where
\begin{equation}\label{8.1}
q(x)=\begin{cases} 1,\quad &\text{if}\quad |x|\le 1\\
|x|^\alpha+|x|^\alpha\cos|x|^\beta,\quad &\text{if}\quad |x|> 1\end{cases}
\end{equation}
under the conditions
\begin{equation}\label{8.2}
\alpha>-2,\qquad \beta>1+\frac{\alpha}{2}.
\end{equation}
Our goal is to use Theorems \ref{thm1.4}, \ref{thm1.8} and \ref{thm7.1} for
finding their analogues in the particular case \eqref{8.1}. For the reader's
convenience, we present the statements proved below as separate theorems
although these ``theorems" are, of course, just examples to the statements
proved above.

\begin{thm}\label{thm8.1}
Suppose that $q(x)$ is of the form \eqref{8.1}. Then for every $x\in R$,
equation \eqref{1.11} has a unique solution $d(x).$ If, in addition, condition
\eqref{8.2} holds, then for all $|x|\gg1,$ we have
\begin{equation}\label{8.3}
d(x)=\frac{1+\delta(x)}{|x|^{\alpha/2}},\qquad
|\delta(x)|\le\frac{c}{|x|^\gamma},\end{equation} where
$\gamma=\min\left\{2,\beta-1-\frac{\alpha}{2}\right\}.
$
\end{thm}
\begin{remark}\label{remark8.1}{\rm
Since the function $q(x)$ in \eqref{8.1} is even, throughout the sequel we will
assume $x\ge0.$ Final results will be written for all $x\in R.$}
\end{remark}

\begin{proof}[Proof of \thmref{thm8.1}]
In the case \eqref{8.1}, relations \eqref{1.2} easily follows from the shape of
the graph of $q(x).$ Then by \lemref{lem2.1}, for every $x\in R $ there exists
a unique positive solution $d(x)$ of equation \eqref{1.11}. To prove formula
\eqref{8.3}, we apply \thmref{thm7.1}. Let $x\gg1,$\ $q_1(x)=x^\alpha,$\
$q_2(x)=x^\alpha\cos x^\beta.$ Then (see \eqref{7.2})
\begin{equation}\label{8.4}
A(x)=[0,2q_1(x)^{-1/2}]=[0,2x^{-\alpha/2}].\end{equation} {}From
\eqref{8.2} for $t\in A(x)$ and $\xi\in[x-t,x+t],$ we get the inequalities
\begin{equation}\label{8.5}
|\xi|\le x+t\le x+2x^{-\alpha/2}=x\left(1+2x^{-1-\frac{\alpha}{2}}\right)\le
3x\quad\text{for}\quad x\gg1,
\end{equation}
\begin{equation}\label{8.6}
|\xi|\ge x+t\ge x-2x^{-\alpha/2}=x\left(1-2x^{-1-\frac{\alpha}{2}}\right)\ge
3^{-1}x\quad\text{for}\quad x\gg1.
\end{equation}
Inequalities \eqref{8.5}--\eqref{8.6} are used below to estimate
$\varkappa_1(x)$ for $x\gg1$ (see \eqref{7.3}):
\begin{align}
\varkappa_1(x)&=\frac{1}{q_1(x)^{3/2}}\sup_{t\in
A(x)}\left|\int_{x-t}^{x+t}q_1''(\xi)d\xi\right|=\frac{1}{x^{3\alpha/2}}\sup_{t\in
A(x)}\left|\int_{x-t}^{x+t}\alpha(\alpha-1)\xi^{\alpha-2}d\xi\right|\nonumber\\
&\le \frac{c}{x^{3\alpha/2}}x^{\alpha-2}\sup_{t\in
A(x)}\left|\int_{x-t}^{x+t}1d\xi\right|=\frac{c(\alpha)}{x^2}.\label{8.7}
\end{align}

Denote $a(x)=\left[x-2x^{\alpha/2},x+2x^{\alpha/2}\right]$ for $x\gg1.$ In the
following estimate for $\varkappa_2(x)$ (see \eqref{7.4}) for $x\gg1,$ we use
relations \eqref{8.5}--\eqref{8.6} and the second mean theorem (\cite[Ch.12,
\S12, no.3]{10}):
\begin{align}
\varkappa_2(x)&=\frac{1}{\sqrt{q_1(x)}}\sup_{t\in
A(x)}\left|\int_{x-t}^{x+t}q_2(\xi)d\xi\right|=
\frac{1}{x^{\alpha/2}}\sup_{t\in
A(x)}\left|\int_{x-t}^{x+t}\xi^{\alpha-\beta+1}\frac{(\beta\xi^{\beta-1}\cos\xi^\beta)d\xi}{\beta}
\right|\nonumber\\
&\le c\frac{x^{\alpha-\beta+1}}{x^{\alpha/2}}\sup_{S_1,S_2\in
a(x)}\left|\int_{S_1}^{S_2}\beta\xi^{\beta-1}\cos\xi^\beta d\xi\right|\le
\frac{c}{x^{\beta-1-\alpha/2}}.\label{8.8}
\end{align}
{}From \eqref{8.2}, \eqref{8.8} and \eqref{8.7}, we get condition \eqref{7.5}.
Now \eqref{8.3} follows from \thmref{thm7.1}.
\end{proof}

\begin{thm}\label{thm8.2}
Suppose that $q(x)$ is of the form \eqref{8.1} and conditions \eqref{8.2} hold.
Then $q(x)\in H.$
\end{thm}

\begin{proof}
Since in this case condition \eqref{1.2} holds (see the proof of
\thmref{thm8.1}), it remains to find a function $k(x)$ satisfying the
requirements of Definition \ref{defn1.2}. Let $m$ be a positive number which
will be chosen later.  Set
\begin{equation}\label{8.9}
k(x)=\begin{cases} 2,\quad & \text{if}\quad |x|\le 2^{\frac{m}{\alpha+2}}\\
|x|^{\frac{\alpha+2}{m}},\quad & \text{if}\quad |x|\ge
2^{\frac{m}{\alpha+2}}\end{cases}\end{equation}

{}From \eqref{8.9} and \eqref{8.2}, it follows that \eqref{1.12} holds. Let us
check \eqref{1.13}. Let $x\gg1$ and $m>2.$ Then from \eqref{8.3} and
\eqref{8.9} it follows that
\begin{equation}\label{8.10}
\frac{k(x)d(x)}{x}\le
c\frac{x^{\frac{\alpha+2}{m}}}{x^{1+\frac{\alpha}{2}}}=cx^{\frac{(\alpha+2)(2-m)}{2m}}\to
0\qquad\text{as}\qquad x\to\iy.
\end{equation}
Therefore from \eqref{8.10} for $t\in[x-k(x)d(x),x+k(x)d(x)]$ and $x\gg1$, we
get
\begin{equation}
\begin{aligned}\label{8.11}
&t\le x+k(x)d(x)=x\left[1+\frac{k(x)d(x)}{x}\right]\le 2x\\
&t\ge x-k(x)d(x)=x\left[1-\frac{k(x)d(x)}{x}\right]\ge
\frac{x}{2}.\end{aligned}
\end{equation}
Inequalities \eqref{1.13} for $x\gg1$ and $t\in[x-k(x)d(x),x+k(x)d(x)]$ follow
from \eqref{8.11}:
\begin{equation}\label{8.12}
\frac{k(t)}{k(x)}=\left(\frac{t}{x}\right)^{\frac{\alpha+2}{m}}\le
2^{\frac{\alpha+2}{m}},\qquad
\frac{k(t)}{k(x)}=\left(\frac{t}{x}\right)^{\frac{\alpha+2}{m}}\ge\left(\frac{1}{2}\right)^
{\frac{\alpha+2}{m}}.
\end{equation}

Estimates \eqref{1.13} for all $x\in R$ can now be derived from \eqref{8.12}
taking into account that the functions under consideration are even and using
\lemref{lem2.6}. Let us check \eqref{1.14}. It is easy to see that in order to
estimate $\Phi(x)$ (see \eqref{1.14}), one can use estimates for $\Phi_1(x),$\
$\Phi_2(x)$ and $\Phi_3(x):$
\begin{align}
\Phi(x)&=k(x)d(x)\sup_{z\in[0,k(x)d(x)]}\left|\int_0^z[q(x+t)-q(x-t)]dt\right|\nonumber\\
&\le
k(x)d(x)\sup_{z\in[0,k(x)d(x)]}\left|\int_0^z[q_1(x+t)-q_1(x-t)]dt\right|\nonumber\\
&\quad +k(x)d(x) \sup_{z\in[0, k(x)d(x)]}\left|\int_x^{x+z}q_2(\xi)d\xi\right|\nonumber\\
&\quad +k(x)d(x) \sup_{z\in[0, k(x)d(x)]}\left|\int_{x-z}^{x}q_2(\xi)d\xi\right|\nonumber\\
&\quad :=\Phi_1(x)+\Phi_2(x)+\Phi_3(x),\qquad x\in R.\label{8.13}
\end{align}

Let $m>6.$ Below in the estimate of $\Phi_1(x)$ for $x\gg1$, we use relations
\eqref{8.11}, \eqref{8.12}, \eqref{8.9} and \eqref{8.3}:
\begin{align}
\Phi_1(x)&=k(x)d(x)\sup_{z\in[0,k(x)d(x)]}\left|\int_0^z\int_{x-t}^{x+t}q_1'(\xi)dt\right|
=k(x)d(x)\sup_{z\in[0,k(x)d(x)]}\left|\int_0^z\int_{x-t}^{x+t}\alpha\xi^{\alpha-1}d\xi
dt\right|\nonumber\\
&\le
ck(x)d(x)x^{\alpha-1}\sup_{z\in[0,k(x)d(x)]}\left|\int_0^z\int_{x-t}^{x+t}d\xi
dt\right|=c(k(x)d(x))^3x^{\alpha-1}\nonumber\\
 &\le c\frac{x^{\frac{3(\alpha+2)}{m}}}{x^{\frac{3\alpha}{2}}}\cdot
 x^{\alpha-1}=c x^{\frac{(\alpha+2)(6-m)}{2m}}\le c.\label{8.14}
 \end{align}

 Since the functions under consideration are even and $\Phi_1(x)$ is continuous
 for $x\in R,$ it is not hard to prove that inequalities \eqref{7.14} (perhaps
 with a bigger constant $c$) hold for all $x\in R.$ Let us now consider
 $\Phi_2(x)$ and $\Phi_3(x).$ We shall prove that these functions are bounded
 for $x\in R$ and $m\gg1;$ since the proof is the same for both functions, below
 we only estimate $\Phi_2(x).$ {}From \eqref{8.2} it follows that there exists
 $m_0$ such that for all $m\ge m_0\ge7$, the following inequalities hold:
 \begin{equation}\label{8.15}
 \beta\ge(\alpha+2)\left(\frac{1}{2}+\frac{1}{m}\right)>\frac{\alpha+2}{2}=1+\frac{\alpha}{2}
 \end{equation}
 \begin{equation}\label{8.16}
 \Rightarrow\quad
 m_0\doe\min_{m\ge7}\left\{m:\frac{(\alpha+2)(m+2)}{2m}\le\beta\right\}.
 \end{equation}

 Denote $b(x)=[x,x+k(x)d(x)].$ Below for $x\gg1$, we use relations
 \eqref{8.11}, \eqref{8.3}, the second mean theorem \cite[Ch.11, \S2, no.3]{9},
 \eqref{8.15} and \eqref{8.16}:
 \begin{align}
  \Phi_2(x)&=k(x)d(x)\sup_{t\in[0,k(x)d(x)]}\left|\int_{x-t}^{x+t}\xi^{\alpha-\beta+1}\frac{
[\beta\xi^{\beta-1}\cos\xi^\beta]}{\beta}d\xi\right|\nonumber\\
 &\le
 ck(x)d(x)x^{\alpha-\beta+1}\sup_{s_1,s_2\in[0,k(x)d(x)]}\left|\int_{s_1}^{s_2}\beta\xi^{
 \beta-1}\cos\xi^\beta d\xi\right|\nonumber\\
 &\le ck(x)d(x)x^{\alpha-\beta+1}\le
 c\frac{x^{\frac{\alpha+2}{m_0}}}{x^{\frac{\alpha}{2}}}x^{\alpha-\beta+1} =
 cx^{\frac{(\alpha+2)(m_0+2)}{2m_0}-\beta}\le c.\label{8.17}
 \end{align}
As in the case of $\Phi_1(x)$ above, one can extend estimate \eqref{8.17} to
the whole axis (perhaps with a bigger constant $c).$ The statement of the
theorem now follows from \eqref{8.13}.
\end{proof}

\begin{cor}\label{cor8.3} Suppose that the function $q(x)$ is defined by
equality \eqref{8.1} under condition \eqref{8.2} and $\rho(x)$ is defined by
equalities \eqref{1.8}. Then for all $|x|\gg1,$ we have the following
asymptotic formula:
\begin{equation}\label{8.18}
\rho(x)=\frac{1+\varepsilon(x)}{2|x|^{\alpha/2}},\qquad |\varepsilon(x)|\le
\frac{c}{|x|^{\gamma_0}}.
\end{equation}
\end{cor}
Here
$\gamma_0=\min\left\{2,\beta-\frac{\alpha}{2}-1,\frac{\alpha+2}{2m_0}\right\}$
and $m_0$ is the number from \eqref{8.16}.

\begin{proof} Formula \eqref{8.18} follows from Theorems \ref{thm8.2},
\ref{thm8.1}, \ref{thm1.4} and the final choice of $k(x)$ made in the course of
the proof of Theorem \ref{thm8.2} and \eqref{1.20}.
\end{proof}

\begin{cor}\label{cor8.4}
Consider the Riccati equation \eqref{1.33} in the case \eqref{8.1} under
condition \eqref{8.2}. The following assertions hold for this equation:
\begin{enumerate}\item[A)] There exists a unique solution $y_1(x)$\ $(y_2(x))$
of equation \eqref{1.32} defined for all $x\in R$ and satisfying the equalities
$$\lim_{x\to-\iy}y_1(x)|x|^{-\alpha/2}=\lim_{x\to\iy}y_1(x)x^{-\alpha/2}=-1$$
$$\left(\lim_{x\to-\iy}y_2(x)|x|^{-\alpha/2}=\lim_{x\to\iy}y_2(x)x^{-\alpha/2}=1\right).$$
\item[B)] Let $y_+(x)$ be a solution of \eqref{1.32} defined on $[c,\iy)$ for
some $c\in R.$ Then $y_+(x)\ne y_1(x)$ for $x\in[c,\iy)$ if and only if
$$\lim_{x\to\iy}y_+(x)x^{-\alpha/2}=1.$$
\item[C)] Let $y_-(x)$ be a solution of \eqref{1.32} defined on $(-\iy,c]$ for
some $c\in R.$ Then $y_-(x)\ne y_2(x)$ for $x\in(-\iy,c]$ if and only if
$$\lim_{x\to-\iy}y_-(x)|x|^{-\alpha/2}=-1.$$
\end{enumerate}
\end{cor}

\begin{proof} This is a consequence of Theorems \ref{thm8.2}, \ref{thm8.1} and
\ref{thm1.8}.
\end{proof}

\end{document}